\documentclass[11pt]{article}
\usepackage{amssymb}

\begin{document}
     
\newtheorem{thm0}{Theorem}[section]
\newtheorem{thm}{Theorem}[subsection]
\newtheorem{thm1}{Theorem}[subsubsection]
\newtheorem{defin0}[thm0]{Definition}
\newtheorem{defin}[thm]{Definition}
\newtheorem{defin1}[thm1]{Definition}
\newtheorem{cor}[thm]{Corollary}
\newtheorem{cor1}[thm1]{Corollary}
\newtheorem{lem}[thm]{Lemma}
\newtheorem{lem1}[thm1]{Lemma}
\newtheorem{prop}[thm]{Proposition}
\newtheorem{rem}[thm]{Remark}
\newtheorem{exam}[thm]{Example}
\newtheorem{exam1}[thm1]{Example}
\newtheorem{sublem}[thm]{Sublemma}
\def\b{{\vrule height6pt width6pt depth2pt}}

\title{K-area, Hofer metric and
geometry of conjugacy classes in Lie groups}

\author{(PRELIMINARY VERSION)\\
\\
Michael Entov\thanks{Partially supported by the 
Nathan, Guta and Robert Gleit\-man Post\-doc\-toral Fellow\-ship
and by the Israel Science Foundation founded by the Israel Academy 
of Sciences and Humanities, grant n. 582/99-1}\\
\\
Department of Pure Mathematics\\
Faculty of Mathematical Sciences\\ 
Weizmann Institute of Science\\
Rehovot 76100\\ 
Israel\\
\\
e-mail: entov@wisdom.weizmann.ac.il\\}

\date{\today}
\maketitle

\begin{abstract}
Given a closed symplectic manifold
$(M,\omega)$
we introduce a certain
quantity associated to a tuple of conjugacy
classes
in the universal cover of the group
${\hbox{\it Ham}}\, (M,\omega)$ 
by means of the Hofer metric on
${\hbox{\it Ham}}\, (M,\omega)$. 
We use pseudo-holomorphic curves
involved in the definition of the multiplicative structure
on the Floer cohomology of a symplectic manifold
$(M,\omega)$
to estimate this quantity 
in terms of actions of some periodic orbits of related
Hamiltonian flows.
As a corollary we get a new way to obtain Agnihotri-Belkale-Woodward
inequalities for
eigenvalues of products of unitary matrices.
As another corollary we get a new proof of the geodesic property 
(with respect to the Hofer metric)
of Hamiltonian
flows generated by certain autonomous Hamiltonians.
Our main technical tool is K-area defined for
Hamiltonian fibrations
over a surface with boundary 
in the spirit of L.Polterovich's work on Hamiltonian fibrations over
$S^2$.

\end{abstract}

%\vfil
%\eject
\tableofcontents
\vfil
\eject

\bigskip
\bigskip
\section{Introduction and overview of the main results}
\label{sec-intro}

\bigskip

\bigskip
\subsection{Definitions and previously known results concerning
$\Delta^G_l$ }
\label{subsect-intro-def-Delta}

With a connected Lie group 
$G$
one can associate an important object
$\Delta^G_l$
defined below.

\bigskip
\begin{defin}
\label{def-Delta}
{\rm
Define
$\Delta^G_l$
as the set of all 
$l{\hbox{\rm -tuples}}$
of conjugacy classes
$({\cal C}_1,\ldots , {\cal C}_l )$
in
$G$
such that
$\varphi_1\cdot\ldots\cdot\varphi_l = Id$
for some
$\varphi_i\in {\cal C}_i$,
$i=1,\ldots, l$.
}
\end{defin}
\smallskip

If
$G = SU (n)$
(see e.g. 
\cite{AW})
the structure of 
$\Delta^G_l$
contains the answer to the following question:
given
$l-1$
matrices 
$A_1,\ldots, A_{l-1}\in SU (n)$,
what can be said about the eigenvalues of the product
$A_1\cdot\ldots\cdot A_{l-1}$
in terms of the eigenvalues of the factors?
In the case
$G = GL (n,{\bf C})$
the structure of
$\Delta^G_l$
is the subject of the so called
{\it
Deligne-Simpson problem
}
and is important in studying possible monodromies of a multi-valued 
solution of a regular Fuchsian system
of differential equations on the Riemannian sphere
--
see e.g.
\cite{Simpson}.

When the Lie group
$G$
is compact, connected and 
finite dimensional it is possible in some cases to get a nice 
description of
$\Delta^G_l$
in geometric terms.
Let
${\mathfrak g}$
be the Lie algebra of
$G$
and let
${\mathfrak t}$
be the Lie algebra of a maximal torus in
$G$.
Choose a fundamental domain 
${\mathfrak U}\subset {\mathfrak t}$
for the adjoint action of
$G$
on
${\mathfrak g}$
so that points in 
${\mathfrak U}$
parameterize all the conjugacy classes in 
$G$:
to each conjugacy class 
${\cal C}$
in
$G$
one associates the unique point
$x\in {\mathfrak U}$
such that
$exp (x)\in {\cal C}$. 
Such a domain
${\mathfrak U}$
is not unique. 
In the case when
$G$
is semi-simple one can choose 
${\mathfrak U}$
as:
\[
{\mathfrak U} = \{ \xi\in {\mathfrak t}_{+}\, |\, \alpha_0 (\xi)\leq 1\}, 
\]
where
${\mathfrak t}$
is the Lie algebra of a maximal torus in
$G$,
${\mathfrak t}_{+}$
is a positive Weyl chamber and
$\alpha_0$
is the highest root.
Thus
$\Delta^G_l$
can be viewed as a subset 
$\Delta^G_l\subset {\mathfrak U }^l\subset {\mathfrak t}^l$.
A theorem of Meinrenken and Woodward
\cite{Me-Wo}
states that if the finite dimensional, compact, connected Lie group
$G$
is also simply connected, then
$\Delta^G_l\subset {\mathfrak t}^l$
is a convex polytope of maximal dimension. 
(This is essentially due to the convexity of
the image of a certain moment map).

If
$G = SU (n)$
the Agnihotri-Belkale-Woodward theorem
(see 
\cite{AW},\cite{Be})
gives a complete description of the convex polytope
$\Delta^G_l$:
the inequalities defining 
$\Delta^G_l$
(we will call them
{\it
ABW inequalities})
are in one-to-one correspondence with non-zero
$l{\hbox{\rm -point}}$
spherical Gromov-Witten invariants of all the complex Grassmannians
$Gr\, (r, n)$,
$1\leq r\leq n-1$.
In other words, the Agnihotri-Belkale-Woodward theorem 
provides a complete list of inequalities 
describing possible eigenvalues of a product of unitary 
matrices in terms of the eigenvalues of the factors.

\bigskip
\subsection{Introducing a pseudo-metric on the group}

In this paper we suggest a new way to study
$\Delta^G_l$.
Namely we equip the tangent bundle of
$G$
with a
{\it
bi-invariant
} 
norm,
i.e. each tangent space is equipped with a norm which varies 
smoothly with the base point.
We will call such a norm
{\it 
Finsler.
}
(Sometimes, given certain
additional assumptions on smoothness and convexity of 
the norm outside of the zero section,
such a structure is called
{\it
an absolutely homogeneous Finsler structure
}
on
$G$
\cite{BCS}).
One can measure lengths of paths in the group 
with respect to a Finsler norm on the tangent bundle and define the 
distance between two points in the group as the
infimum of lengths of paths connecting them. 
This defines a 
{\it
Finsler
bi-invariant pseudo-metric
}
$\rho$
on
$G$, 
i.e. it is a bi-invariant symmetric function on
$G\times G$
which satisfies the triangular inequality but may vanish not only on
the diagonal but also outside of it.
If 
$\rho\, (x,y)\neq 0$
as long as
$x\neq y$
then the pseudo-metric
$\rho$
is a genuine 
{\it
Finsler metric 
}
on
$G$.
If
$G$
is finite-dimensional then such a Finsler norm on
$T_\ast G$ 
always defines a genuine
Finsler metric on the group. 
In the infinite-dimensional case this may not be the case.

An important example of a group with a bi-invariant Finsler norm
is provided a certain
$C^0{\hbox{\rm -norm}}$ 
on the tangent bundle of the infinite-dimensional Lie group
${\hbox{\it Ham}}\, (M,\omega)$
defining the famous Hofer metric on the group
--
see the definitions in
Section~\ref{subsect-intro-main-result}.
(More precisely, we will mostly consider the universal cover 
${\widetilde{\hbox{\it Ham}}}\, (M,\omega)$
of
${\hbox{\it Ham}}\, (M,\omega)$
equipped with the pullback of the norm on the tangent bundle of
${\hbox{\it Ham}}\, (M,\omega)$
under the covering map).

\bigskip
\begin{defin}
\label{def-upsilon}
{\rm
Given a bi-invariant Finsler pseudo-metric
$\rho$
on
$G$
and a tuple
${\cal C} = ({\cal C}_1,\ldots , {\cal C}_l)$
of conjugacy classes in
$G$
define
\[
\Upsilon_l ({\cal C}) = 
\inf_{\varphi_i\in 
{\cal C}_i } \rho\, (Id, \varphi_1\cdot\ldots\cdot\varphi_l).
\]
}
\end{defin}
\bigskip

The quantity
$\Upsilon_l$
has the following elementary properties.

\medskip
\noindent
1) If
$\Upsilon_l\, ({\cal C}) = 0$
then
${\cal C}\in \Delta^G_l$
and if
$G$
is finite-dimensional then 
$\Delta^G_l = 
\{ \, {\cal C} = 
({\cal C}_1,\ldots , {\cal C}_l ) \, | 
\, \Upsilon_l\, ({\cal C}) = 0 \, \}$.

\smallskip
\noindent
2) 
$\Upsilon_l \, ({\cal C}_1,\ldots , {\cal C}_l ) = 
\inf_{{\cal C}_{l+1}} \rho\, (Id, {\cal C}_{l+1})$,
where the infimum is taken over all 
${\cal C}_{l+1}$
such that
$({\cal C}_1,\ldots , {\cal C}_l, {\cal C}_{l+1})\in\Delta^G_{l+1}$.

\smallskip
\noindent
3) $\Upsilon_l$
does not depend on the ordering of classes in
a tuple -- it is actually a function on 
${\it sets}$ 
of conjugacy classes.

\smallskip
\noindent
4) 
$\Upsilon_2\, ({\cal C}_1, {\cal C}_2) = 
\rho\, ({\cal C}_1^{-1}, {\cal C}_2 )$, 
where the right-hand side denotes the distance between the two closed
sets in the pseudo-metric space.

\smallskip
\noindent
5)
$\Upsilon_l\, ({\cal C}_1,\ldots ,{\cal C}_l ) =
\Upsilon_l\, ({\cal C}_1^{-1},\ldots , {\cal C}_l^{-1})$.

\smallskip
\noindent
6) Triangular inequality:
$\Upsilon_l\, ({\cal C}_1,\ldots , {\cal C}_l )\geq 
\rho\, (Id, {\cal C}_l) -
\rho\, (Id, {\cal C}_{l-1}) - \ldots - 
\rho\, (Id, {\cal C}_1).
$
\bigskip

Observe that the property 1) is true no matter which Finsler pseudo-metric
$\rho$
on
$G$
we take. Thus a good choice of the bi-invariant Finsler pseudo-metric 
may help us in understanding the structure of
$\Delta^G_l$.
Such a good choice will come from the pullback
of the Finsler norm on the tangent bundle of 
${\hbox{\it Ham}}\, (M,\omega)$
defining the Hofer metric under a homomorphism
$G\to {\hbox{\it Ham}}\, (M,\omega)$
representing a Hamiltonian action of
$G$
on a closed symplectic manifold
$(M,\omega)$.
As an example consider the case of
$G = SU (n)$ 
acting on
$({\bf C}P^{n-1}, \omega)$, 
where
$\omega$
is the standard (Fubini-Studi) symplectic form.
The action is Hamiltonian and thus it defines 
a homomorphism
$SU (n)\to {\hbox{\it Ham}}\, ({\bf C}P^{n-1},\omega)$
with a the kernel
${\bf Z}_n$.
The pullback of the Finsler norm on the tangent bundle of
${\hbox{\it Ham}}\, ({\bf C}P^{n-1},\omega)$
is a bi-invariant Finsler norm on
$SU (n)$ 
defined by the operator norm
$\| A\| = \sup_{\| x\| =1}  \| A(x)\|$
on the Lie algebra of
$SU (n)$.

If
$G$
is an infinite-dimensional Lie group it is not clear how to describe
$\Delta^G_{l+1}$
as a geometric object.
Nevertheless, as it will be shown below, in the case when
$G = {\widetilde{\hbox{\it Ham}}}\, (M,\omega)$
one can use a certain estimate on
$\Upsilon_l\, ({\cal C}_1,\ldots , {\cal C}_l )$
to get a generalization of the ABW inequalities. 
The same sort of estimate also leads to some interesting 
applications concerning the geodesics on
${\hbox{\it Ham}}\, (M,\omega)$
with respect to the Hofer metric.
The estimate on
$\Upsilon_l\, ({\cal C}_1,\ldots , {\cal C}_l )$
will be given in terms of actions of some periodic orbits of Hamiltonian
flows related to
${\cal C}_1,\ldots , {\cal C}_l$.

\vfil
\eject

\bigskip
\subsection{The main result}
\label{subsect-intro-main-result}

Before stating the result we briefly recall basic definitions.

\bigskip
\noindent
{\bf Hamiltonian functions, vector fields and symplectomorphisms.}

Let
$(M,\omega)$
be a closed symplectic manifold.
A function
$h: S^1\times M\to {\bf R}$
(called
{\it
Hamiltonian function
} 
or simply
{\it
Hamiltonian
})
defines a (time-dependent)
{\it
Hamiltonian vector field 
}
$X_h$
by the formula
\begin{equation}
\label{eqn-def-ham-vect-field}
d h^t (\cdot) = \omega (X_{h^t}, \cdot),
\end{equation}
where 
$h^t  = h (t,\cdot)$,
and the formula holds pointwise on
$M$ 
for every
$t\in S^1$.

A 
{\it
Hamiltonian symplectomorphism
}
of
$(M,\omega)$
is a diffeomorphism of
$M$
which can be represented as the time-1 map of the 
flow of a (time-dependent) Hamiltonian vector field
$X_h$.
In such a case we say that the Hamiltonian symplectomorphism is
{\it
generated
}
by the Hamiltonian
$h$.
\bigskip

\bigskip
\noindent
{\bf Normalization condition.}
All the Hamiltonians
$h:S^1\times M\to {\bf R}$
in the paper will be assumed to be normalized so that
$\int_M h^t \omega^n = 0$,
for any
$t\in S^1$.
Such a Hamiltonian will be called 
{\it
a normalized Hamiltonian.
}
\bigskip

\bigskip
\noindent
{\bf The group 
${\hbox{\it Ham}}\, (M,\omega)$
and the Hofer metric on it.}

Hamiltonian symplectomorphisms of 
$(M,\omega)$
form a group
${\hbox{\it Ham}}\, (M,\omega)$
which can be viewed as an infinite-dimensional Lie group:
it can be equipped with the structure of 
an infinite-dimensional manifold so that the group product 
and taking the inverse of an element become smooth operations
\cite{Ra-Schm}.
The Lie algebra 
${\cal H}$
of 
${\hbox{\it Ham}}\, (M,\omega)$
can be identified with the (Poisson) Lie algebra 
of all functions on 
$M$
with the zero mean value.
The norm
$\| h \| = \max_M | h | $
on 
${\cal H}$
defines a bi-invariant Finsler norm on the tangent bundle
of the group
${\hbox{\it Ham}}\, (M,\omega)$.
It is a deep result of symplectic topology 
(see
\cite{Ho}, 
\cite{Pol-Hm},
\cite{L-McD-Hm})
that the Finsler norm on the tangent bundle of the group
leads to a genuine bi-invariant Finsler metric 
$\rho$
on the group itself. This metric is called 
{\it
Hofer metric
}
(the original metric introduced by Hofer was actually defined by the norm
$\| h\| = \max_M h - \min_M h$
and is equivalent to the metric we use).
The Hofer metric on
${\hbox{\it Ham}}\, (M,\omega)$
lifts to a bi-invariant Finsler pseudo-metric on the universal cover
${\widetilde{\hbox{\it Ham}}}\, (M,\omega)$
of the group
${\hbox{\it Ham}}\, (M,\omega)$.
The Finsler pseudo-metric on
${\widetilde{\hbox{\it Ham}}}\, (M,\omega)$
defines the corresponding function
$\Upsilon_l$
on 
$l{\hbox{\rm -tuples}}$
of conjugacy classes in
${\widetilde{\hbox{\it Ham}}}\, (M,\omega)$.
The Hamiltonian flow of
$X_h$
over a period of time from 0 to 1, 
viewed as a path in
${\hbox{\it Ham}}\, (M,\omega)$
starting at
$Id$,
determines an element in 
${\widetilde{\hbox{\it Ham}}}\, (M,\omega)$.
We will denote the conjugacy class of such an element in
${\widetilde{\hbox{\it Ham}}}\, (M,\omega)$
by
${\cal C}_h$.
\bigskip

\bigskip
\noindent
{\bf Action functional.}

Consider the set of pairs 
$(\gamma, f)$,
where
$\gamma: S^1\to M$
is a contractible curve and
$f: D^2\to M$, 
is a disk spanning the curve
$\gamma$, 
i.e. 
${\left. f\right|}_{\partial D^2} = \gamma$.
Introduce an equivalence relation on the set of such pairs as follows:
two pairs
$(\gamma, f)$
and
$(\gamma, f^\prime )$
are called 
{\it
equivalent
}
if the connected sum
$f\sharp f^\prime$
represents a torsion class in
$H_2 (M, {\bf Z})$.
Denote by
${\cal P} (h)$
the space of equivalence classes 
${\hat\gamma} = [\gamma, f ]$
of pairs
$(\gamma, f)$,
where
$\gamma: S^1\to M$
is a contractible time-1 periodic trajectory of the Hamiltonian flow of
$h$.
The theorems proving the Arnold's conjecture 
\cite{Fu-Ono}, \cite{Liu-Ti}
imply that any Hamiltonian symplectomorphism of
$M$
has a closed contractible periodic orbit and therefore
${\cal P} (h)$
is always non-empty.

Define 
$\Pi$
as the group of spherical homology classes of
$M$, 
i.e. as the image of the
Hurewicz homomorphism 
$\pi_2 (M)\to H_2 (M, {\bf Z})/{\rm Tors}$.
The group
$\Pi$
acts on 
${\cal P} (h)$
by the formula:
\[
A : [\gamma , f]\mapsto [\gamma, A\sharp f].
\]
For an element
${\hat\gamma} = [\gamma, f]\in {\cal P} (h)$
define its 
{\it
action:
}
\[
{\cal A}_h ({\hat\gamma}) = 
- \int_D f^\ast \omega - \int_0^1 h (t, \gamma (t)) dt.
\]

Given 
$l$ 
Hamiltonians 
$H = (H_1,\ldots ,H_l)$
on
$M$
we denote by
${\cal P} (H)$
the set of equivalence classes
${\hat\gamma} = [{\hat{\gamma}}_1, \ldots ,{\hat{\gamma}}_l]$,
where
${\hat{\gamma}}_i = [\gamma_i, f_i] \in {\cal P}(H_i)$
and the equivalence relation is given by
\[
[{\hat{\gamma}}_1, \ldots ,{\hat{\gamma}}_l ]\sim 
[A_1 \sharp {\hat{\gamma}}_1, \ldots , A_l \sharp {\hat{\gamma}}_l ],
\]
whenever
$A_i\in \Pi$
and
$A_1+\ldots +A_l$
is a torsion class.
The group
$\Pi$
acts on
${\cal P} (H)$:
\[
A: [{\hat{\gamma}}_1, \ldots ,{\hat{\gamma}}_l ]\mapsto
[A\sharp {\hat{\gamma}}_1, {\hat{\gamma}}_2,\ldots ,{\hat{\gamma}}_l ].
\]
Set
${\cal A}_H ({\hat\gamma}) 
= {\cal A}_{H_1} ({\hat\gamma}_1)+\ldots + {\cal A}_{H_l} ({\hat\gamma}_l)$
for
${\hat\gamma}\in {\cal P} (H)$.
\bigskip

\bigskip
\noindent
{\bf Moduli spaces
${\cal M} (\hat{\gamma}, H, \tilde{J})$.
}

Let
$\Sigma$
be a complex Riemann surface of genus 0 with
$l$
cylindrical ends
$\Sigma_i\cong [0,+\infty )\times S^1$,
$i=1,\ldots ,l$.
Consider the trivial fibration
$\Sigma\times M\to\Sigma$
and let
$pr_M$
denote the natural projection of
$\Sigma\times M$
on
$M$.
Let
$H = (H_1, \ldots , H_l)$
be some (time-dependent) Hamiltonians and let
${\hat\gamma}_i = [\gamma_i, f_i ]\in {\cal P} (H_i)$,
$i=1,\ldots ,l$.
Set
$\hat{\gamma} = (\hat{\gamma}_1,\ldots , \hat{\gamma}_l)$.
We will be interested in a certain class 
${\cal T} (H)$
of almost complex structures on
$\Sigma\times M$
depending on
$H$ --
see
Section~\ref{sect-psh-curves-gw-invariants}
for the definition.
Given an almost complex structure
$\tilde{J}\in {\cal T} (H)$
we consider the 
$\tilde{J}{\hbox{\rm -holomorphic}}$ 
sections 
$u: \Sigma\to \Sigma\times M$
such that the ends of the surface
$pr_M\circ u (\Sigma)\subset M$
converge at infinity to the periodic orbits 
$\gamma_1,\ldots ,\gamma_l$
and such that the curve 
$pr_M\circ u (\Sigma)\subset M$
capped off with the discs
$f_1 (D^2), \ldots, f_l (D^2)$
forms a closed surface representing a torsion integral homology class in
$M$.
The space of such pseudo-holomorphic sections 
$u$
will be denoted by
${\cal M} (\hat{\gamma}, H, \tilde{J})$.

Given (time-dependent) normalized Hamiltonians 
$H= (H_1,\ldots, H_l)$
on
$M$
we will introduce a certain number
${\hbox{\it size}}\, (H)$
(see 
Definition~\ref{def-sigma-weak-coupl})
and consider a certain family of subsets 
$\{ {\cal T}_\tau (H)\} \subset {\cal T} (H)$,
$0< \tau < {\hbox{\it size}}\, (H)$ 
-- see
Definition~\ref{def-cal-T-tau}.
We will say that 
$\hat{\gamma}\in {\cal P} (H)$
is 
{\it
durable
}
if there exists a sequence
$\{ \tau_k \}\nearrow {\hbox{\it size}}\, (H)$
such that all the spaces
${\cal M} (\hat{\gamma}, H, \tilde{J}_{\tau_k})$
are non-empty.

Now we state our main result.

%\vfil
%\eject

\bigskip
\begin{thm}
\label{thm-exist-psh-curve-implies-estimate}
Let
$H = (H_1,\ldots, H_l)$
be some (time-dependent) normalized Hamiltonians on
$M$
and let
${\cal C}_H = ({\cal C}_{H_1}, \ldots , {\cal C}_{H_l} )$ 
be the corresponding conjugacy classes in
${\widetilde{\hbox{\it Ham}}}\, (M,\omega)$.
Then for any durable 
$\hat{\gamma}\in {\cal P} (H)$
one has 
\[
\Upsilon_l\, ({\cal C}_H ) \geq 
{\cal A}_H (\hat{\gamma}).
\]
\end{thm}
\medskip

The inequalities from 
Theorem~\ref{thm-exist-psh-curve-implies-estimate}
will be called
{\it
action inequalities.
}
Each inequality boils down to the fact that the integral of a 
certain symplectic form
$\Omega_\tau$
on
$\Sigma\times M$
over a
$\tilde{J}{\hbox{\rm -holomorphic}}$
curve is non-negative as soon as
$\tilde{J}$
is compatible with the symplectic form
$\Omega_\tau$
(i.e.
$\Omega_\tau (\cdot, \tilde{J}\cdot )$
defines a Riemannian metric on
$\Sigma\times M$).

The fact that certain moduli spaces 
${\cal M} (\hat{\gamma}, H, \tilde{J})$
are non-empty can be checked in some cases by considering the
structure of the ``pair-of-pants'' product
in the Floer cohomology on the level of Floer cochains --
the zero-dimensional spaces 
${\cal M} (\hat{\gamma}, H, \tilde{J})$
are instrumental in defining such a product
(see
Section~\ref{subsect-floer-quantum-multiplication}).

\bigskip
\begin{cor}
\label{cor-abw-ineqs-ham-usual-form}
Suppose that 
$\hat{\gamma}\in {\cal P} (H)$
is durable and
$\Upsilon_l\, ({\cal C}_H ) = 0$.
Then
${\cal A}_H (\hat{\gamma}) \leq 0$.
\end{cor}
\smallskip

This is a generalization of the ABW inequalities for the case of the Lie
group
$G = {\widetilde{\hbox{\it Ham}}}\, (M,\omega)$.
However, unlike in
\cite{AW}, \cite{Be},
we do not know whether the inequalities from
Corollary~\ref{cor-abw-ineqs-ham-usual-form}
provide a complete list of inequalities defining
$\Delta_l^G$.

\bigskip
\subsection{An overview of the applications}

In
Section~\ref{subsect-Grassmannians-ABW-ineqs}
we will demonstrate how the ABW inequalities for
$SU (n)$
can be recovered from action inequalities if one considers the natural
group actions of 
$SU (n)$
on all the complex Grassmannians
$Gr\, (r, n)$,
$1\leq r\leq n-1$,
and if the elements of
$SU (n)$
are viewed as Hamiltonian symplectomorphisms of the Grassmannians.
This approach also indicates
what should be the generalization of the ABW inequalities
for other Lie groups
(see 
Section~\ref{subsect-Grassmannians-ABW-ineqs}).

We will also apply
Theorem~\ref{thm-exist-psh-curve-implies-estimate}
to the study of geodesics in the group
${\hbox{\it Ham}}\, (M,\omega)$,
extending previous results of Lalonde and McDuff
(see Section~\ref{subsect-time-indep-hams-applications}).

\bigskip
\subsection{K-area}

The key observation which provides the connection between the function
$\Upsilon_l$
coming from the Hofer geometry on
${\widetilde{\hbox{\it Ham}}}\, (M,\omega)$
and the pseudo-holomorphic curves is that 
$\Upsilon_l ({\cal C}_1, \ldots ,{\cal C}_l )$
can be interpreted in terms of a certain
{\it
K-area.
}
Roughly speaking, K-area is the inverse of 
a quantity obtained by fixing a class
of connections on some fibration and taking infimum
of a
$C^0{\hbox{\rm -norm}}$ 
of the curvature tensors of connections in the fixed class (assuming
that we have some prefixed metrics on the base and on the fiber used 
to measure the norm of a curvature tensor).
The notion of K-area and its applications to symplectic topology
were first introduced by M.Gromov 
in his seminal paper 
\cite{Gro}.
Other remarkable applications of K-area 
to symplectic topology were discovered later by L.Polterovich
(see
\cite{Pol1}, \cite{Pol2}, \cite{Pol3}).
He studied Hamiltonian fibrations over
$S^2$
and, in particular, found a close connection between the following objects:

\medskip
\noindent
$\bullet$
the K-area of a Hamiltonian fibration over
$S^2$ 
(where the bi-invariant 
Hofer metric on the group of Hamiltonian symplectomorphisms
of the fiber is used for the measurements in the definition of 
the K-area);

\smallskip
\noindent
$\bullet$
the Hofer length of the clutching loop for the Hamiltonian
fibration (i.e. the loop in the group of Hamiltonian symplectomorphisms
of the fiber used to construct the fibration over
$S^2$
from the trivial fibrations over the two hemispheres).
\medskip

\noindent
L.Polterovich also found a way to use K-area as a tool
to extract an estimate on the Hofer length of 
the clutching loop from the fact that a symplectic form
on the total space of the fibration integrates non-negatively 
over a pseudo-holomorphic section of the fibration 
(see 
\cite{Pol1}).

In this paper we extend Polterovich's methods and results to
fibrations 
over an oriented surface 
$\Sigma$
with boundary.
We will be especially interested in the case when the genus of 
$\Sigma$
is zero and the fibration is Hamiltonian, i.e. its typical fiber is a closed
symplectic manifold
$(M,\omega)$
and the structural group is
$G = {\hbox{\it Ham}}\, (M,\omega)$.
(For another interesting development of the Polterovich's ideas in the 
similar spirit but in a different setup see 
\cite{Akv-Sal}). 

A
$G{\hbox{\rm -fibration}}$
with a connected fiber over a surface
$\Sigma$
as above is always trivial. Therefore in order to have 
an interesting quantity one needs to define K-area using the connections 
whose holonomies along the boundary components of
$\Sigma$
lie in some fixed conjugacy classes 
$({\cal C}_1,\ldots , {\cal C}_l)$.
Thus K-area becomes a function that associates a number to each tuple 
$({\cal C}_1,\ldots , {\cal C}_l)$
of conjugacy classes in 
$G$, 
and as it will be shown in the paper, in the case when the genus of
$\Sigma$
is zero, this function is actually the 
inverse of
$\Upsilon_l$.
(In particular, if
$G$
is finite-dimensional then the tuples from
$\Delta^G_l$
are exactly all possible tuples of conjugacy classes of holonomies
of flat connections
over
$\partial \Sigma$
on a trivial fibration
$G\times\Sigma\to\Sigma$).

\bigskip
\noindent
{\bf 
Acknowledgments.
}
I am deeply grateful to L.Polterovich who suggested the problem to me and
generously shared with me his ideas on the subject. 
I am greatly indebted to D.McDuff and the anonymous 
referee who carefully read an 
earlier version of 
the manuscript, pointing out errors and making valuable suggestions.
I also thank 
J.Bernstein, 
Y.Eliashberg, A.Givental, F.Lalonde, D.Salamon for useful discussions,  
M.Schwarz 
for an important comment and for sending me his PhD-thesis and A.Klyachko for 
pointing out to me the crucial references. This work was started during 
the author's stay at Tel Aviv University. I thank it for the hospitality.

\bigskip

\vfil
\eject

\bigskip
\bigskip
\section{Applications of
Theorem~\ref{thm-exist-psh-curve-implies-estimate}
}
\label{sect-main-appl-Ham}

\bigskip
\subsection{Preliminaries}
\label{sect-appl-of-main-thm-prelim}

Before stating the applications of
Theorem~\ref{thm-exist-psh-curve-implies-estimate}
we recall the basic facts about quantum and Morse cohomology
and define a special class of time-independent Hamiltonians,
called
{\it
slow Hamiltonians.
}

\bigskip
\subsubsection{Strongly semi-positive symplectic manifolds}
\label{subsect-semi-positive}

We say that an almost complex structure
$J$
on
$(M,\omega)$
is
{\it 
compatible 
}
with
$\omega$
if
$\omega (\cdot, J\cdot)$
is a Riemannian metric on
$M$.
Almost complex structures compatible with
$\omega$
form a contractible space.

For the technical reasons coming from the theory of pseudo-holomorphic curves
we assume from now on that 
$(M^{2n},\omega)$
is 
{\it
strongly semi-positive,
}
although in view of the recent developments
(see 
\cite{Fu-Ono},
\cite{Liu-Ti},
\cite{Lu})
it is likely that this assumption can be removed.
Namely, a closed symplectic manifold
$(M^{2n},\omega)$
is called
{\it
strongly semi-positive,
} 
if for every
$A\in \pi_2 (M)$
one has
\begin{equation}
\label{eqn-str-semipos}
2-n\leq c_1 (A)< 0 \Longrightarrow \omega (A)\leq 0,
\end{equation}
where
$c_1$
is the first Chern class of the tangent bundle 
$T_\ast M$
equipped with an almost complex structure
$J$
In 
\cite{Se}
such manifolds were said to satisfy the 
``Assumption
$W^{+}$''.

Recall that 
$M$
is called 
{\it 
semi-positive 
}
\cite{PSS}
or
{\it
weakly monotone
}
\cite{Ho-Sa}
if in 
(\ref{eqn-str-semipos})
one replaces the inequality
$2-n\leq c_1 (A)\leq 0$
by
$3-n\leq c_1 (A)\leq 0$
thus weakening the condition.
Strong semi-positivity of
$M$
means that for 
$m=0,1,2,3$
a generic 
$m{\hbox{\rm -parametric}}$
family of 
$\omega{\hbox{\rm -compatible}}$
almost complex structures on
$(M,\omega )$
does not contain an almost complex structure
$J$
for which there exist
$J{\hbox{\rm -holomorphic}}$
spheres with a negative Chern number. 
In the case of usual (not strong) semi-positivity this holds only for
$m=0,1,2$
and since we would like to consider 2-parametric families of
almost complex structures and 1-parametric deformations of such families
(as in 
\cite{Se})
we ask the manifold to be strongly semi-positive.

The class of strongly 
semi-positive symplectic manifolds includes closed symplectic
surfaces, symplectic tori, complex projective spaces, complex
Grassmannians, complex flag manifolds and other interesting objects.
\bigskip

\bigskip
\subsubsection{Quantum cohomology}
\label{subsect-quantum-basic-facts}

Let
$(M,\omega)$,
be a closed connected strongly semi-positive symplectic manifold.
We will briefly recall the necessary basic facts about quantum
cohomology of
$(M,\omega)$.

First, we define the coefficient ring for the quantum cohomology.
Namely, from the group
$\Pi$
one builds an appropriate Novikov graded ring
$\Lambda_\omega$
over 
${\bf Q}$:
an element of
$\Lambda_\omega$
is a formal sum
\[
\lambda = \sum_{A\in\Pi} \lambda_A e^{2\pi i A}
\]
with rational coefficients
$\lambda_A\in {\bf Q}$
which satisfies the condition
\[
\sharp \{ A\in\Pi\, |\, \lambda_A\neq 0, \omega (A) \leq c\} <\infty
\]
for every
$c > 0$.
The natural multiplication makes
$\Lambda_\omega$
a ring. This ring has a natural grading defined by
${\rm deg}\, (e^{2\pi i A}) = 2 c_1 (A)$.
The subring of elements of degree zero will be denoted by
$\Lambda_0\subset \Lambda_\omega$.

The
{\it
quantum cohomology group
}
of
$M$
is defined as a graded tensor product
$QH^\ast (M) = H^\ast (M) \otimes \Lambda_\omega$,
where
$H^\ast (M)$
denotes the quotient
$H^\ast (M,{\bf Z})/Tors$.

Given homology classes
$A\in\Pi$,
and
$\alpha_j\in H_{i_j} (M)$,
$j=1,\ldots , l$,
satisfying the condition
\[
i_1+\ldots +i_l = 2n (l-1) - 2 c_1 (A)
\]
one can define the
{\it
Gromov-Witten number
}
$(\alpha_1,\ldots,\alpha_l)_A$
(see 
\cite{Wi},
\cite{Ru-Ti},
\cite{Ru-Ti-1}).
The Gromov-Witten numbers define the ring structure on the quantum
cohomology group.
Namely, let
$a_1,\ldots, a_l\in H^\ast (M)$
and let
$PD (a_1), \ldots,$ 
$PD (a_{l-1})$
be the corresponding Poincar{\'e}-dual homology classes.
The multiplication, or 
{\it
the quantum product
}
on the quantum cohomology group is defined by the formula:
\[
a_1\ast \ldots\ast a_{l-1} = \sum_{A\in\Pi} (a_1\ast\ldots\ast a_{l-1})_A 
e^{2\pi i A},
\]
where the class
$(a_1\ast\ldots\ast a_{l-1})_A \in H^{(2n-i_1) +\ldots 
+ (2n-i_{l-1}) - 2c_1 (A)} (M)$
has to satisfy the condition
\[
\langle (a_1 \ast\ldots\ast a_{l-1})_A, \alpha_l\rangle = 
(PD (a_1), \ldots, PD (a_{l-1}),\alpha_l)_A
\]
for any homology class
$\alpha_l\in H_\ast (M)$.
It was proved in 
\cite{Ru-Ti}
that the quantum product is associative.
The class 
$(a_1\ast\ldots\ast a_l)_0$
represents the usual cup-product:
$(a_1\ast\ldots\ast a_l)_0 = a_1\cup\ldots\cup a_l$.
The cohomology class 
${\bf 1}\in H^0 (M)$
Poincar{\'e}-dual to the fundamental class
$[M]$
is the unit element in
$QH^\ast (M)$.

\bigskip
\subsubsection{Morse homology and cohomology}

Here we briefly the definitions of Morse homology and cohomology
(see e.g.
\cite{Sch-book}
for details).
Fix a Riemannian metric
on
$M^{2n}$.
Suppose that
$h$
is a Morse function on
$M^{2n}$,
and moreover, it is a Morse-Smale function with respect to
the metric.
In such a case one can define the
{\it
Morse chain complex
}
$C_\ast (h)$
of the function
$h$
over the graded coefficient ring 
$\Lambda_\omega$
as a free graded module over
$\Lambda_\omega$
generated by the critical points of
$h$
graded by their Morse indices.
The differential in such a complex is 
defined by means of counting downward
gradient trajectories of
$h$
with respect to the metric
that connect the critical points of neighboring indices
(see 
\cite{Sch-book}
for details).
Similarly one can define the dual cochain complex 
$C^\ast (h)$
of
$C_\ast (h)$.
The homology of the chain complexes 
$C_\ast (h)$
and 
$C^\ast (h)$
are called, respectively,
{\it
the Morse homology and cohomology 
}
of
$h$.
The Morse homology (resp. cohomology) of a Morse function 
is canonically isomorphic to the singular homology (resp. cohomology) of
$M$.
The tautological identification
$C^k (-h) \cong C_{2n-k} (h)$
leads to the Poincar{\'e} isomorphism between the Morse cohomology
of
$-h$
and the Morse homology of
$h$.

\bigskip
\subsubsection{Critical points of functions and multiplicative 
identities in quantum cohomology}
\label{subsubsect-crit-pts-slow-hams-multipl-ident-quant-cohom}

\bigskip
\begin{defin}[cf. \cite{Pol-book}]
\label{def-hom-ess-cr-pt}
{\rm
Fix a Riemannian metric
on the closed symplectic manifold
$M^{2n}$.
Let
$h$
be a Morse-Smale function on 
$M$
with respect to the metric.
Let 
$Crit\, (h)$
denote the set of all critical points of
$h$,
or the set of all generators of the Morse complex
$C_\ast (h)$
over
$\Lambda_\omega$.
whose differential is denoted by
$\partial$.
Identify the Morse homology of 
$h$
with the singular homology of
$M$.
We say that a critical point 
$y_i$
of
$h$
{\it
homologically essential
}
for a homology class 
$\alpha\in H_\ast (M,\Lambda_\omega)$
if 
$\alpha$,
viewed as a class in
$H_\ast (C_\ast (h) )$,
does not belong to the image of
$i_{{\cal K}_\ast}: H_\ast ({\cal K}_\ast )\to H_\ast (C_\ast (h) )$
for any subcomplex
${\cal K}_\ast\subset {\rm Span}\, (Crit\, (h)\setminus \{ y_i\})$.
Equivalently, one can say that 
$y_i$
is homologically essential for a homology class
$\alpha$
if it enters with a non-zero coefficient into any chain in
$C_\ast (h)$
representing 
$\alpha$.

}
\end{defin}
\smallskip

Similar questions concerning
the necessity of critical points in a chain representing a
homology class in the context of generating functions 
were studied by C.Viterbo in 
\cite{Vit}.

\bigskip
\begin{exam}[cf. \cite{Pol-book}]
\label{exam-max-hom-ess}
{\rm
Suppose that
$y$
is a unique point of 
{\it 
global 
}
maximum for a function
$h$
which is Morse-Smale with respect to some Riemannian metric on
$M$.
Let
$(C_\ast (h) , \partial)$
be the Morse complex of 
$h$.
Then we claim that
$y$
is homologically essential for the fundamental class
$[M]\in H_\ast (C_\ast (h) )$.
Indeed, the subspace
${\rm Span}\, (Crit\, (h)\setminus \{ y \} )\subset C_\ast (h)$
is
$\partial{\hbox{\rm -invariant}}$
and the chain complex
$({\rm Span}\, (Crit\, (h)\setminus \{ y\}), \partial )$
is nothing else but the Morse complex for the function
$h$
on an
{\it 
open
}
manifold
$M\setminus y$.
But
$H_{2n} (M) = \Lambda_\omega$
and
$H_{2n} (M\setminus y) =  0$
which proves the claim.

}
\end{exam}
\smallskip

Let
$H_1,\ldots, H_l$
be Morse Hamiltonians on
$(M,\omega )$.
Let
$A\in\Pi$
and let
$z^i$,
$i= 1,\ldots, l$,
be a critical point of
$H_i$.
We say that
${\hat{\gamma}} = [{\hat{\gamma}}_1 ,\ldots ,  {\hat{\gamma}}_l ]\in 
{\cal P} (H)$
is
{\it
associated with
}
$z^1,\ldots, z^l, A$
if
${\hat{\gamma}}_i\in {\cal P} (H_i)$,
$i=1,\ldots ,l-1$,
is formed by the pair of constant maps to
$z^i$
and
${\hat{\gamma}}_l\in {\cal P} (H_i)$
is formed by a constant map
$S^1\to z^l$
to 
$z^l$
and a smooth 2-sphere attached to
$z^l$
which realizing the spherical homology class
$A$.
\smallskip

\bigskip
\noindent
{\bf Convention: identification of Morse and singular (co)homology.}
\smallskip

\noindent
From now on, given Morse Hamiltonians
$H_1,\ldots , H_l$
we will always identify the singular homology and cohomology of
$M$
with, respectively, the Morse homology and cohomology of
$-H_i$,
$i=1,\ldots ,l$.
\bigskip

\bigskip
\begin{defin}
\label{def-involved}
{\rm
Suppose that for some cohomology classes 
$c_1,\ldots, c_{l-1}\in H^\ast (M, {\bf Q})$
one has
\begin{equation}
\label{eqn-identity-to-be-involved-in}
c_1\ast\ldots\ast c_{l-1} = \sum_{B\in\Pi} c_B e^{2\pi i B},
\end{equation}
where
$c_B\in H^\ast (M, {\bf Q})$.
Let 
$z^i$,
be a critical point of
$H_i$,
$i=1,\ldots , l$,
and let
$A\in\Pi$.

We say that
$\hat{\gamma}$
associated with
$z^1,\ldots ,z^l, A$
is
{\it
involved
}
in the identity
(\ref{eqn-identity-to-be-involved-in})
if the following conditions hold:

\medskip
\noindent
$\bullet$
The point
$z^l$,
viewed as a critical point of  
$-H_l$,
is homologically essential for the rational (singular) homology class
Poincar{\'e}-dual to
$c_A$,
viewed, under our identification convention, 
as a homology class of the Morse complex
$C_\ast (-H_l)$.

\medskip
\noindent
$\bullet$
For each
$i=1,\ldots, l-1$
there exists a basis
$\alpha_1^i, \ldots , \alpha_N^i$,
$N= {\hbox{\rm dim}}\, H_\ast (M, {\bf Q})$,
of
$H_\ast (M, {\bf Q})$
over
${\bf Q}$
such that exactly one basic element 
$\alpha_{j(i)}^i$
satisfies the following two conditions:

\smallskip
(A) 
$z_i$, 
viewed as a critical point of
$-H_i$,
is essential for
$\alpha_{j (i)}^i$,
viewed, under our identification convention, 
as a homology class of the Morse complex
$C_\ast (-H_i)$.

\smallskip
(B)
$ c_i (\alpha_{j (i)}^i)\neq 0$.
\smallskip

Such a basis
$\{ \alpha_j^i\}$
of
$H_\ast (M, {\bf Q})$
will be called
$i{\hbox{\it -friendly}}$.

}
\end{defin}
\smallskip

\bigskip
\begin{exam}
\label{exam-tuples-of-maxs-min-involved-in-multipl-ident}
{\rm
Consider the element
${\bf 1}\in QH^0 (M)$
and the identity
\begin{equation}
\label{eqn-product-of-l-1-identities}
{\bf 1}\ast\ldots\ast {\bf 1} = {\bf 1}
\end{equation}
with 
$l-1$
factors in the left-hand side.
Suppose that
$z_i$,
$i =1,\ldots, l-1$,
is a unique point of global maximum of
$H_i$
and
$z_l$
is a unique point of global minimum of
$H_l$.
Let
${\hat{\gamma}}_i$,
$i=1,\ldots, l$,
be formed by a pair of constant maps into
$z_i$.
Then,
as one easily sees from
Example~\ref{exam-max-hom-ess},
${\hat{\gamma}}$
is involved in the identity
(\ref{eqn-product-of-l-1-identities}).
In this particular case any basis in homology is 
$i{\hbox{\rm -friendly}}$
for any
$i=1,\ldots ,l-1$, 
since 
it has to include exactly one generator from
$H_0 (M, {\bf Q})\cong {\bf Q}$.

}
\end{exam}
\smallskip

\bigskip
\begin{exam}
\label{exam-tuples-of-crit-pts-involved-in-multipl-ident}
{\rm
Suppose that
$H_1,\ldots, H_l$
are perfect Morse functions on
$M$,
i.e. the differentials in their
Morse chain complexes over the integers, 
and hence in the Morse complexes
$C_\ast ( - H_i )$,
$C^\ast ( - H_i )$,
$i=1,\ldots ,l$,
over
$\Lambda_\omega$,
are zero. 
Denote the critical points of
$- H_i$,
$i=1,\ldots ,l$,
by
$z^i_1, \ldots, z^i_N$,
$N= {\hbox{\rm dim}}\, H_\ast (M, {\bf Q})$.
Identify
$QH_\ast (M)$,
as before with the Morse homology of
$C_\ast ( - H_i )$,
$i=1,\ldots ,l$
with the coefficients in
$\Lambda_\omega$.

For a fixed 
$i$
the points 
$z^i_j$,
$j=1, \ldots , N$, 
viewed as the homology classes in
$H_\ast (C_\ast ( - H_i ))$
form a basis of
$H_\ast (C_\ast ( - H_i ))\cong QH_\ast (M)$
over
$\Lambda_\omega$.
Denote the dual basis in the cohomology
$H_\ast (C^\ast ( - H_i ))\cong QH^\ast (M)$
by
$Z^i_j$,
$j= 1, \ldots , N$,
i.e.
$Z^i_j$
takes value 1 on
$z^i_j$
and zero on any other critical point of
$-H_i$.

The cohomology classes 
$Y^i_j$,
$j= 1, \ldots , N$,
Poincar{\'e}-dual to
$z^i_j$,
form another basis in 
$H_\ast (C^\ast ( - H_i ))\cong QH^\ast (M)$,
and the homology classes
$y^i_j$,
$j= 1, \ldots , N$,
Poincar{\'e}-dual to
$Z^i_j$,
form another basis in the homology
$H_\ast (C_\ast ( - H_i ))\cong QH_\ast (M)$
dual to the basis
$Y^i_j$,
$j= 1, \ldots , N$, 
in the cohomology.
In such a case we will say that 
the 
{\it
homology class
$y^i_j$
is
Poincar{\'e}-dual 
to
$z^i_j$.
}

Now let
$A\in\Pi$
and let
$z^i_{j (i)}$,
$i=1,\ldots, l$,
be a critical point of
$- H_i$
of Morse index
$m (i)$,
so that
\[
\sum_{i=1}^l (2n -m(i)) - 2c_1 (A) = 2n.
\]
Let
${\hat{\gamma}}_i$,
$i=1,\ldots, l-1$,
be formed by a pair of constant maps into
$z^i_{j (i)}$
and let
${\hat{\gamma}}_l\in {\cal P} (H_i)$
be formed by a constant map
$S^1\to z^l_{j (l)}$
and a smooth 2-sphere attached to
$z^l_{j (l)}$
that realizes the homology class 
$A$
(here
$z^i_{j (i)}$,
$i=1,\ldots ,l$,
is viewed as a critical point of
$H_i$).

Write the quantum product of the cohomology classes
$Z^1_{j (1)},\ldots , Z^{l-1}_{j (l-1)}$
as
\begin{equation}
\label{eqn-identity-Z}
Z^1_{j (1)}\ast\ldots\ast Z^{l-1}_{j (l-1)} = 
\sum_{B\in \Pi} \lambda_B e^{2\pi iB} = 
\sum_{B\in\Pi, 1\leq j\leq N} \lambda_{B,j} Y^l_j e^{2\pi iB},
\end{equation}
where each
$\lambda_B\in H^\ast (M,\Lambda_\omega)$
is decomposed along the basis
$\{ Y^l_j\}$,
$j= 1, \ldots , N$,
with the coefficients
$\lambda_{B,j}\in {\bf Z}$.
According to the definition of quantum multiplication, the coefficient
$\lambda_{A, j (l)}$
equals to the Gromov-Witten number
$(y^1_{j (1)},\ldots , y^l_{j(l)} )_A$.

Then
${\hat{\gamma}} = ({\hat{\gamma}}_1,\ldots, {\hat{\gamma}}_l )$
is involved in the identity
(\ref{eqn-identity-Z})
if and only if the Gromov-Witten number
$(y^1_{j (1)},\ldots , y^l_{j (l)} )_A$,
or the coefficient
$\lambda_{A,j}$
at the term
$Y^l_{j(l)} e^{2\pi i A}$
in
(\ref{eqn-identity-Z})
is non-zero.

In this case
for each
$i=1,\ldots ,l-1$
the homology classes
$z^i_j$,
$j= 1, \ldots , N$,
form an
$i{\hbox{\rm -friendly}}$ 
basis of
$H_\ast (M, {\bf Q})$.

Observe that this example can be immediately generalized to the situation when
not all the differentials in the Morse complexes
$C_\ast ( - H_i )$,
$i=1,\ldots ,l$,
are zero but only the ones surrounding the level
$m (i)$:
\[
\partial_{m (i) +1}: C_{m (i) +1} (-H_i ) \to C_{m (i)} (-H_i ) 
\]
and
\[
\partial_{m (i)}: C_{m (i)} (- H_i )\to C_{m (i) -1} (-H_i ). 
\]

}
\end{exam}
\medskip

\subsubsection{Slow Hamiltonians}
\label{sect-slow-hams-def}

\bigskip
\begin{defin}
\label{def-slow-hams}
{\rm
A Hamiltonian 
$h$
on a closed symplectic manifold 
$(M,\omega)$
is called
{\it
slow
}
if it satisfies the following conditions:

\medskip
\noindent
$(A)$
$h$
is time-independent;

\smallskip
\noindent
$(B)$
the Hamiltonian flow of
$h$
has only constant contractible periodic 
trajectories of period less or equal than
$1$;

\smallskip
\noindent
$(C)$
the Hessian of
$h$
at any of its critical points does not have an imaginary eigenvalue
$i\lambda$
with
$\lambda\geq 2\pi$.

}
\end{defin}
\smallskip

A (local) result of Siegel and Moser 
\cite{Sieg-Mo}
(also see
\cite{McD-Slim})
shows that a generic Hamiltonian satisfying the conditions
(A) and (B) of 
Definition~\ref{def-slow-hams}
also satisfies the condition
(C).

As an example of a slow Hamiltonian one can pick any
sufficiently
$C^2{\hbox{\rm -small}}$
function on
$M$
(see e.g.
\cite{HZ}).

\bigskip
\subsection{Time-independent Hamiltonians: the main result}
\label{subsect-time-indep-hams-main-result}

The following result is essentially based on 
Theorem~\ref{thm-exist-psh-curve-implies-estimate}
and on the fact that for a
slow Morse Hamiltonian the Floer complex can be identified with the Morse
complex and the pair-of-pants product on the dual Morse complex
descends to the quantum product on the Morse cohomology -- see
Section~\ref{subsect-floer-quantum-multiplication}.

\bigskip
\begin{thm}
\label{thm-main-time-indep-hams}
Let
$(M,\omega)$
be strongly semi-positive.
Let
$H = (H_1,\ldots, H_l)$
be normalized slow Morse Hamiltonians.
Let
$z^i$,
$i=1,\ldots ,l$,
be some critical points of
$H_i$
and let
$A\in\Pi$
so that
$\hat{\gamma} \in {\cal P} (H)$,
associated with
$z^1,\ldots, z^l, A$,
is involved in some identity
(\ref{eqn-identity-to-be-involved-in})
in the quantum cohomology of
$(M,\omega)$.

\noindent
Then
\[
\Upsilon_l\, ({\cal C}_H ) \geq 
{\cal A}_H (\hat{\gamma}) = -H_1 (z^1) -\ldots - H_l (z^l) - \omega (A).
\]

\end{thm}
\smallskip

Using
Example~\ref{exam-tuples-of-crit-pts-involved-in-multipl-ident}
one immediately gets from 
Theorem~\ref{thm-main-time-indep-hams}
the following 

\bigskip
\begin{cor}
\label{cor-thm-main-time-indep-hams-perfect-hams}
Let
$H = (H_1,\ldots, H_l)$
be normalized slow Hamiltonians which are perfect Morse functions on a
strongly semi-positive
$(M,\omega)$.
Let
$A\in\Pi$
and let
$z^i$,
$i=1,\ldots ,l$,
be a critical point of
$-H_i$
of Morse index
$m (i)$,
viewed under the identification
$H_\ast (C_\ast ( - H_i ))\cong QH_\ast (M)$
as an integral homology class of
$M$
of degree
$m (i)$,
so that
\[
m (1) +\ldots + m(l) = 2n (l-1) - 2c_1 (A).
\]
Let
$y^i\in H_{2n- m(i)} (M)$,
$i=1,\ldots ,l$,
be the homology class Poincar{\'e}-dual to
$z^i$
and suppose that
$(y^1,\ldots , y^l)_A\neq 0$.

\noindent
Then
\[
\Upsilon_l\, ({\cal C}_H ) \geq 
{\cal A}_H (\hat{\gamma}) = -H_1 (z^1) -\ldots - H_l (z^l) - \omega (A).
\]
\end{cor}
\smallskip

\bigskip
\subsection{Application: geodesics in 
${\hbox{\it Ham}}\,  (M,\omega)$.
}
\label{subsect-time-indep-hams-applications}

Assume as before that
$(M,\omega)$
is a closed strongly semi-positive symplectic manifold.

Suppose that 
$H_1, H_2: M\to {\bf R}$
are normalized slow Hamiltonians that generate, respectively, 
Hamiltonian symplectomorphisms
$\varphi_{H_1}$
and
$\varphi_{H_2}$
\break
whose conjugacy classes in
${\hbox{\it Ham}}\,  (M,\omega)$
are denoted by
$[\varphi_{H_1}]$
and
$[\varphi_{H_2}]$.
Denote by
${\cal C}_{H_1}$
and
${\cal C}_{H_2}$
the conjugacy classes in the universal cover of
${\hbox{\it Ham}}\,  (M,\omega)$
corresponding, respectively, to the Hamiltonian flows of
$H_1$
and
$H_2$
over the interval of time from 0 to 1.
Let 
$\{a\}$
be a homotopy class of paths in
${\hbox{\it Ham}}\,  (M,\omega)$
connecting
$[\varphi_{H_1}]$
and
$[\varphi_{H_2}]$
whose lifts to the universal cover of
${\hbox{\it Ham}}\,  (M,\omega)$
connect
${\cal C}_{H_1}$
and
${\cal C}_{H_2}$.
Denote by
$\rho_{\{ a \}} ([\varphi_{H_1}] , [\varphi_{H_2}])$
the infimum of lengths of paths connecting 
$[\varphi_{H_1}]$
and
$[\varphi_{H_2}]$
from the homotopy class
$\{a\}$.

\bigskip
\begin{thm}
\label{thm-two-small-hams}
In the notation as above one has
\[
\rho_{\{ a \}  } ([\varphi_{H_1}] , [\varphi_{H_2}])\geq
\max (| \max_M H_1 - \max_M H_2 |, | \min_M H_1 - \min_M H_2 |).
\]
\end{thm}
\smallskip

Given two Hamiltonian flows 
$\{ f_t\}, \{ g_t \}$,
generated, respectively, by Ham\-il\-tonians 
$F, G$,
the composition
$\{ h_t = f_t g_t\}$
is a Hamiltonian flow generated by the Hamiltonian
$H (t, x) = F(t, x) + G(t, f_t^{-1} (x) )$. 
Also given a Hamiltonian symplectomorphism
$\varphi_H$
generated by a time-independent Hamiltonian
$H: M\to {\bf R}$
one has that
$\phi^{-1}\circ\varphi_H\circ\phi $
can be generated by the Hamiltonian function
$H\circ\phi$
for any
$\phi\in {\hbox{\it Ham}}\,  (M,\omega)$.
This immediately imposes an estimate from above on 
$\rho_{\{ a \}  } ([\varphi_{H_1}] , [\varphi_{H_2}])$:
\[
\rho_{\{ a \}  } ([\varphi_{H_1}] , [\varphi_{H_2}])\leq 
\inf_{\phi\in {\hbox{\it Ham}}\,  (M,\omega)} \| H_1 - H_2\circ \phi\, \|,
\]
where
$\|\cdot\|$
is the norm 
$\| h\| = \max_M | h|$
on the Lie algebra of
${\hbox{\it Ham}}\,  (M,\omega)$
defining the Hofer metric on the group.

Let
$H$
be a normalized slow Hamiltonian.
Considering 
$H_1 = \tau_1 H , H_2 = \tau_2 H$
for any
$0< \tau_1 <\tau_2\leq 1$,
and applying
Theorem~\ref{thm-two-small-hams}
to the Hamiltonians
$H_1, H_2$
and then to the Hamiltonians 
$-H_1, - H_2$
generating the symplectomorphisms
$\varphi_{-H_1} = \varphi^{-1}_{H_1}$,
$\varphi_{-H_2} = \varphi^{-1}_{H_2}$,
one readily obtains the following corollary.

\bigskip
\begin{cor}
\label{cor-auton-ham-gives-geod}
For any
$0\leq\tau_1,\tau_2\leq 1$
the path
$\{\varphi_H^\tau\}$,
$\tau\in [\tau_1,\tau_2]$,
is globally length-minimizing in its homotopy class in
${\hbox{\it Ham}}\,  (M,\omega)$.
\end{cor}
\smallskip

This result has been previously proved by F. Lalonde and D.McDuff
for symplectic manifolds which are either 2-dimensional or weakly exact
(see
\cite{L-McD},
Theorem 5.4).
An extension of this result to the case of a general closed
symplectic manifold has also been the subject of a recent independent
work by D.McDuff and J.Slimowitz
\cite{McD-Slim}
(also see
\cite{Sli}).

Let us now state a result similar to
Theorem~\ref{thm-two-small-hams}
for the case when
$l>2$.

\bigskip
\begin{thm}
\label{thm-many-small-hams}
Suppose that 
$H= (H_1,\ldots, H_l)$
are normalized slow Ham\-il\-tonians and let
${\cal C}_H = ({\cal C}_{H_1},\ldots , {\cal C}_{H_l} )$
be the corresponding conjugacy classes in
${\widetilde{\hbox{\it Ham}}}\, (M,\omega)$.
Then
\begin{equation}
\label{ineq-upsilon-ham-case-many-hams}
\Upsilon_l ({\cal C}_H ) \geq 
-\max_M H_1 -\ldots - \max_M H_{l-1} -\min_M H_l .
\end{equation}
\end{thm}
\smallskip

As an easy application of 
Theorem~\ref{thm-two-small-hams}
consider the case when
$M = S^2$
with the standard symplectic form
$\omega$
normalized so that
$\int_{S^2} \omega = 1$.
Let
$\varphi_1$
and
$\varphi_2$
be the linear rotations of
$S^2$
around the same axis by angles
$0\leq\xi_1\leq\xi_2\leq\pi$
respectively.
Clearly,
$\varphi_1, \varphi_2\in {\hbox{\it Ham}}\, (S^2,\omega)$.
In this case we can actually compute
$\rho_{\{ a\}  } ([\varphi_1], [\varphi_2 ])$
for all possible homotopy classes 
$\{ a\}$:
there are four such homotopy classes.
Indeed, according to a result of Gromov 
\cite{Gro-pshc}
the subgroup
$SO (3)\subset {\hbox{\it Ham}}\, (S^2,\omega)$
is a deformational retract of
${\hbox{\it Ham}}\, (S^2, \omega )$.
Hence
$\pi_1 ({\hbox{\it Ham}}\, (S^2, \omega )) = \pi_1 (SO (3) ) = {\bf Z}_2$
and the path in
$SO (3)$
defining the full
$2\pi$
twist of
$S^2$
around some axis represents the generator of 
${\bf Z}_2$.
This shows that there only two homotopy classes of paths connecting
each of the elements
$\varphi_1$
and
$\varphi_2$
with the identity in
${\hbox{\it Ham}} (S^2,\omega)$
(they are represented by ``clockwise'' and ``counterclockwise''
rotations around the oriented axis).
Therefore we have four possible homotopy classes
$\{  a\}$.
All the Hamiltonians can be easily computed:
$\varphi_i$
is generated by an autunomous Hamiltonian
$H_i : S^2 \to {\bf R}$
which has only two critical points with critical values
$\max H_i = \zeta_i /2$
and
$\min H_i = - \zeta_i /2$.

Applying 
Theorem~\ref{thm-two-small-hams}
we recover the result of
F. Lalonde and 
\break
D.McDuff
(see
\cite{L-McD},
Corollary 1.10)
concerning geodesics in
${\hbox{\it Ham}}\,  (S^2,\omega)$
(also see
\cite{Pol1}
for further developments).

\bigskip
\begin{cor}[\cite{L-McD}]
\label{cor-geod-sphere}
The distance 
$\rho\, (\varphi_1, \varphi_2)$
between 
$\varphi_1$ 
and 
$\varphi_2$
in the group 
${\hbox{\it Ham}}\,  (S^2)$
equals
$\displaystyle{\frac{\zeta_2 -\zeta_1}{2}}$.
\end{cor}
\smallskip

As another application of 
Theorem~\ref{thm-exist-psh-curve-implies-estimate}
we give an estimate from below on
$\Upsilon_3$
in the case of the symplectic two-dimensional torus
$M={\bf T}^2 = {\bf R}^2 /{\bf Z}^2$.
This estimate is different from the one in 
Theorem~\ref{thm-many-small-hams}.
Observe that in the case of
$M={\bf T}^2$
there is no need to pass to the universal cover of
${\hbox{\it Ham}}\,  ({\bf T}^2,\omega)$
since this group is already simply connected
(see e.g. 
\cite{Pol-book}).

Let
$H_1, H_2, H_3: {\bf T}^2\to {\bf R}$
be normalized slow Hamiltonians which are perfect Morse functions.
Denote by
$[\varphi_{H_1} ], [\varphi_{H_2} ], [\varphi_{H_3}]$
the conjugacy classes of the Hamiltonian symplectomorphisms generated,
respectively, by 
$H_1, H_2$, 
$H_3$.

Denote by
$z_3$
the point of global maximum of
$H_3$
and denote by
$z_1 , z_2$
the critical points of index 1 of, respectively,
$H_1$
and 
$H_2$
so that 
$z_1$
and
$z_2$
correspond to the homology classes 
$\alpha_1, \alpha_2$
that generate
$H_1 ({\bf T}^2 )$.

\bigskip
\begin{thm}
\label{thm-3-slow-hams-two-dim-torus}
In the notation as above
\[
\Upsilon_3\, ([\varphi_{H_1} ], [\varphi_{H_2} ], [\varphi_{H_3}])
\geq 
- H_1 (z_1) - H_2 (z_2) - H_3 (z_3).
\]
\end{thm}
\smallskip

\bigskip
\subsection{Application: Grassmannians, ABW inequalities}
\label{subsect-Grassmannians-ABW-ineqs}

In this section we discuss the relation between the ABW and the action
inequalities. 

Consider the Grassmannian 
$Gr\, (r,n)$
of complex
$r{\hbox{\rm -planes}}$
in
${\bf C}^n$.
The Grassmannian
can be viewed as the result of symplectic reduction
for the Hamiltonian action of
$U(r)$
(by multiplication from the right) on the space of complex
$r\times n$
matrices.
The symplectic structure 
$\omega$
on
$Gr\, (r,n)$
that we choose is constructed as the one induced
on the symplectic reduction by the standard 
symplectic structure on the space
${\bf C}^{nr}$
of complex
$r\times n$
matrices.
The natural complex structure on 
$Gr\, (r, n)$
is
compatible with
$\omega$
and we assume that 
$\omega$ 
is normalized so that the cohomology class
$[\omega ]$
takes value 1 on the generator of
$H_2 (Gr\, (r, n), {\bf Z})= \Pi\cong {\bf Z}$
realized by a complex sphere.
One quickly checks that
$(Gr\, (r,n),\omega)$
is a strongly semi-positive symplectic manifold.

If one fixes a complete flag 
$\{0\} = F_0\subset\ldots \ldots\subset F_n = {\bf C}^n$
then to each subset
$I = \{ i_1,\ldots ,i_r\}\subset \{ 1,\ldots , n \}$
corresponds a Schubert variety 
\[
\{ W_I\in Gr\, (r,n)\ |
\ {\rm dim}\, (W_I\cap F_{i_j})\geq j,\, j=1,\ldots ,r \}
\]
representing an integer homology class
$\sigma_I$
which does not actually depend on the choice of the flag.
Given homology classes
$\sigma_{I_1},\ldots, \sigma_{I_l}\in H_\ast (Gr\, (r,n))$
and a class
$d\in {\bf Z}\cong H_2 (Gr\, (r,n), {\bf Z})$
one can define the corresponding Gromov-Witten invariant
$(\sigma_{I_1},\ldots, \sigma_{I_l})_d$.
The Gromov-Witten invariants mentioned above determine the multiplicative
structure of the quantum cohomology ring of
$(Gr\, (r,n),\omega )$
(see
\cite{Bertr}, 
\cite{Bertr-CiF-Fult}, 
\cite{Sieb-Ti},
\cite{Wit-grasm}).

The group
$SU (n)$
naturally acts on
$Gr\, (r,n)$.
This action is Hamiltonian and it produces 
the homomorphism
$SU (n)\to {\hbox{\it Ham}}\, (Gr\, (r,n),\omega)$
with the 
{\it
finite
}
kernel
${\bf Z}_n$.
Since the kernel is finite the homomorphism is a covering map from
$G$
to its image 
which is an embedded submanifold of
${\hbox{\it Ham}}\, (M,\omega)$.
This allows to pull back the Finsler
norm on 
$T_\ast {\hbox{\it Ham}}\, (Gr\, (r,n),\omega)$
(defining the Hofer metric on the group
${\hbox{\it Ham}}\, (Gr\, (r,n),\omega)$)
to an induced Finsler norm on
$T_\ast SU (n)$.
This Finsler norm leads to a Finsler metric on
$SU (n)$
and by means of this metric
one defines the function
$\Upsilon_l$
on 
$l{\hbox{\rm -tuples}}$
of conjugacy classes in
$SU (n)$.
Also consider the pullback of the same Finsler norm from
$T_\ast {\hbox{\it Ham}}\, (Gr\, (r,n),\omega)$
to
$T_\ast {\widetilde{\hbox{\it Ham}}}\, (Gr\, (r,n),\omega) $
and consider the function
$\Upsilon^{(r)}_l$
on 
$l{\hbox{\rm -tuples}}$
of conjugacy classes in
${\widetilde{\hbox{\it Ham}}}\, (Gr\, (r,n),\omega) $
defined by means of the corresponding Finsler pseudo-metric
on the group.

Choose the fundamental domain 
${\mathfrak U}$
in the Cartan subalgebra of
$su (n)$
to be defined by the equations:
\[
\alpha_1 +\ldots +\alpha_n =0,\ 
\alpha_1\geq\ldots\geq\alpha_n\geq \alpha_1 -1,
\]
so that
${\cal C}_\alpha$,
$\alpha = (\alpha_1 ,\ldots ,\alpha_n )\in {\mathfrak U}$,
is the conjugacy class of an element of
$SU (n)$
with the eigenvalues
$e^{2\pi i \alpha_1},\ldots ,e^{2\pi i\alpha_n}$.

Given
$\alpha = (\alpha_1,\ldots ,\alpha_n)\in {\mathfrak U}$
consider the 1-parametric subgroup 
\break
$diag (e^{2\pi i\alpha_1 t},\ldots ,e^{2\pi i\alpha_n t})
\subset SU (n)$
of diagonal matrices.
The action of this 1-parametric subgroup on
$Gr\, (r,n)$
produces a Hamiltonian flow
$\{ \varphi_{H_\alpha}^t \}$
for some time-independent Hamiltonian
$H_\alpha: Gr (r,n)\to {\bf R}$
with mean value zero.
Thus
$H_\alpha$
is normalized.
Denote by
${\cal C}^{(r)}_\alpha := {\cal C}_{H_\alpha}$
the conjugacy class in
${\widetilde{\hbox{\it Ham}}}\, (Gr\, (r,n),\omega) $
containing the element defined by the flow.

One can check that for
$\alpha = (\alpha_1,\ldots ,\alpha_n)\in {\mathfrak U}$
such that
$\alpha_i \neq \alpha_j$,
$1\leq i,j\leq n$,
the Hamiltonian
$H_\alpha$
is a Morse function. Such a Hamiltonian 
$H_\alpha$
has two crucial properties
(see
Section~\ref{sect-pf-thm-abw-ineqs}):

\bigskip
\noindent
1) 
$H_\alpha$
is a slow Hamiltonian. This is basically due to the
the fact that
${\mathfrak U}$
lies in the domain of injectivity of the exponential
map. 

\medskip
\noindent
2)
$H_\alpha$
is a perfect Morse function.
Indeed, the critical points of
$H_\alpha$
are in one-to-one correspondence with the invariant complex
$r{\hbox{\rm -dimensional}}$
subspaces of the matrix
$diag (e^{2\pi i\alpha_1 t},\ldots ,e^{2\pi i\alpha_n t})$
or, equivalently, with the subsets
$I= \{ i_1,\ldots ,i_r\}\subset \{ 1,\ldots , n \}$, 
or, equivalently, with the free generators
$\sigma_I$
of the homology group of
$Gr\, (r,n)$.
\medskip

These properties of
$H_\alpha$
will allow us to identify the Floer chain complex corresponding to 
$H_\alpha$
and the Morse chain complex
corresponding to 
$- H_\alpha$
so that the quantum product of cohomology classes of 
$(Gr\, (r,n), \omega)$
gets identified with the Floer product of the same classes realized by
the Morse cohomology classes of
$-H_\alpha$
(see
Section~\ref{subsect-floer-quantum-multiplication}).
\bigskip

Now let
$\zeta^i = (\zeta_1^i,\ldots ,\zeta_n^i)\in {\mathfrak U}$,
$i=1,\ldots ,l$,
$\zeta = (\zeta^1 ,\ldots , \zeta^l)\in {\mathfrak U}^l$
and let
${\cal C}_\zeta = ({\cal C}_{\zeta^1},\ldots ,{\cal C}_{\zeta^l} )$
be the corresponding conjugacy classes in
$SU (n)$.
For the conjugacy classes 
${\cal C}_\zeta = ({\cal C}_{\zeta^1},\ldots ,{\cal C}_{\zeta^l} )$
in
$SU (n)$
denote by
${\cal C}^{(r)}_\zeta = 
({\cal C}^{(r)}_{\zeta^1},\ldots , {\cal C}^{(r)}_{\zeta^l})$
the corresponding conjugacy classes in 
${\widetilde{\hbox{\it Ham}}}\, (Gr\, (r,n),\omega) $.

\bigskip
\bigskip
\begin{thm}
\label{thm-abw-ineqs}
Suppose 
$(\sigma_{I_1},\ldots ,\sigma_{I_l})_d\neq 0$
in 
$Gr\, (r,n)$.
Then
\[
\Upsilon_l ({\cal C}_\zeta ) \geq
\Upsilon^{(r)}_l ({\cal C}^{(r)}_\zeta )\geq
\sum_{j=1}^l \sum_{i\in I_j} \zeta_i^j - d.
\]

\end{thm}
\smallskip

As a corollary we obtain the original ABW inequalities for
$SU (n)$.

\bigskip
\begin{cor}[cf. \cite{AW},\cite{Be}]
\label{cor-abw-ineqs}
Suppose
$\Upsilon_l ({\cal C}_\zeta ) = 0$.
Then for any
$r$,
$1\leq r\leq n-1$,
and any
$I_j$,
$|I_j | = r$,
$j=1,\ldots ,l$,
such that
$(\sigma_{I_1},\ldots ,\sigma_{I_r})_d\neq 0$
in
$Gr\, (r,n)$
one has
$\sum_{j=1}^l \sum_{i\in I_j} \zeta_i^j \leq d$.
\end{cor}
\smallskip

In fact, according to
\cite{AW},\cite{Be}
the converse is also true:
if for any
$r$,
$1\leq r\leq n-1$,
and any
$I_j$,
$| I_j | = r$,
$j=1,\ldots ,l$,
such that
$(\sigma_{I_1},\ldots ,\sigma_{I_r})_d\neq 0$
in
$Gr\, (r,n)$
one has
$\sum_{j=1}^l \sum_{i\in I_j} \zeta_i^j \leq d$
then
$\Upsilon_l ({\cal C}_\zeta ) = 0$.
(To see it observe that all the ABW inequalities for all
$l$
can be deduced from the inequalities for
$l=2$,
in which case the statement can be checked easily).

\bigskip
\begin{exam}
\label{exam-Delta-Upsilon-SU2}
{\rm
Consider the case
$G = SU (2)$,
$l=3$.
The conjugacy class of a matrix from 
$SU (2)$
with eigenvalues
$e^{\pm 2\pi i\zeta}$,
$0\leq\zeta\leq 1/2$,
is completely determined by the real number
$\zeta$.
Then
$\Delta^G_3$
is polytope of maximal dimension which lies inside the cube
$[0,1/2]\times [0,1/2]\times [0,1/2]$
in
${\bf R}^3$.
The inequalities defining
$\Delta^G_3$
and
$\Upsilon_3$
can be computed in this case directly by elementary methods 
(cf. \cite{Je-We})
and the result, of course,
matches 
Theorem~\ref{thm-abw-ineqs}
and
Corollary~\ref{cor-abw-ineqs}.
The polytope
$\Delta^G_3$
is a tetrahedron defined by the equations
\[
\zeta^1 +\zeta^2 +\zeta^3\leq 1,
\] 
\[
\zeta^1\leq \zeta^2 +\zeta^3,
\]
\[
\zeta^2\leq \zeta^1 +\zeta^3,
\]
\[
\zeta^3\leq \zeta^1 +\zeta^2,
\]
corresponding to its four faces.
The function
$\Upsilon_3$
on a triple of conjugacy classes (viewed each as a point in
$[0,1/2]$)
is given by the formula:
\[
\Upsilon_3 (\zeta^1, \zeta^2, \zeta^3) =
\max \{ 0, \zeta^3 - \min (\zeta^1 +\zeta^2, 1-\zeta^1 -\zeta^2)\}.
\]
Considering this formula for
$(\zeta^1, \zeta^2, \zeta^3)$
close to zero one gets that the inequalities in 
Theorem~\ref{thm-abw-ineqs}
may turn into equalities.

}
\end{exam}
\smallskip

The proofs of the results above 
indicate how one can possibly describe the convex polytope
$\Delta^G_l$
for a compact semi-simple connected and 
simply-connected Lie group
$G$
other than
$SU (n)$.
Let us briefly sketch how this can 
be done. Given such a Lie group
$G$
with the Lie algebra
${\mathfrak g}$
one should consider all its compact K{\"a}hler homogeneous spaces
(these are in one-to-one correspondence with subsets of the set
of simple roots of
$G$ --
see 
\cite{Serre}, \cite{Wang}).
Since
$G$
is semi-simple, one has 
$[{\mathfrak g} , {\mathfrak g} ] = {\mathfrak g}$
and
$H^2 ({\mathfrak g}) = 0$.
Therefore the natural action of
$G$
on any such homogeneous space 
$(M,\omega )$
is Hamiltonian. (As F.Lalonde has pointed out to me, this is true even if 
$G$ is not semi-simple: using the flux homomorphism one can easily deduce it from the fact that $G$ is just compact and simply-connected). The kernel of this action is the center
$G$
which is finite because $G$ is semi-simple.
Since 
$[{\mathfrak g} , {\mathfrak g} ] = {\mathfrak g}$
any element of
$G$
acts on
$(M,\omega )$
as a Hamiltonian symplectomorphism 
generated by a Hamiltonian
(see 
Definition~\ref{def-slow-hams})
that can be represented as the Poisson bracket 
of some other Hamiltonians and therefore is normalized. 
Thus the Hamiltonian action of 
$G$
on
$(M,\omega )$
induces an inclusion of
${\mathfrak g}$
into the Lie algebra of functions on
$M$
with mean value zero.
Pick a fundamental domain
${\mathfrak U}\subset {\mathfrak t}\subset {\mathfrak g}$
containing zero and lying in the domain of injectivity of the exponential
map.
Now given an element
$v\in {\mathfrak U}$
consider the Hamiltonian flow on
$(M,\omega )$
induced by the action of the 1-parametric subgroup
$\{ e^{\tau v}\}_{ 0\leq\tau\leq 1}$.
For a generic
$v$
this Hamiltonian flow is generated by a normalized
Morse Hamiltonian function
$H_v$.
As in the case of
$SU (n)$
the function
$H_v$
satisfies two crucial properties:

\medskip
\noindent
1) 
$H_v$
is a slow Hamiltonian. As in the case of
$SU (n)$
this is based on the fact that
${\mathfrak U}$
lies in the domain of injectivity of the exponential
map. 

\smallskip
\noindent
2)
$H_v$
is a perfect Morse function. This can be checked using
the cell decomposition for
$M$
which is the analogue of the Schubert cell 
decomposition for complex Grassmannians
(see
\cite{Borel}, 
\cite{Chev}).
\medskip

The actions of the constant periodic trajectories of the
Hamiltonian flow generated by
$H_v$
depend only on the conjugacy class of
$e^{\tau v}$
in
$G$
and completely determine that class.
Then to any non-zero Gromov-Witten number 
$(\alpha_1,\ldots, \alpha_l )_A$
for generators 
$\alpha_1,\ldots, \alpha_l$ 
from the basis of 
$H_\ast (M)$ 
given by the cell decomposition
one associates the corresponding action inequality 
for tuples of conjugacy classes in the universal cover of
${\hbox{\it Ham}}\, (M,\omega)$.
This would in turn provide a generalized ABW inequality that
would have to be satisfied by any tuple belonging to
$\Delta^G_l$.
Considering all the non-zero Gromov-Witten numbers for all 
compact K{\"a}hler homogeneous spaces 
$(M,\omega)$
of
$G$
one should get a complete set of inequalities defining the
convex polytope
$\Delta^G_l$.

Recently, in
\cite{Tel-Wood},
C.Teleman and C.Woodward, using the method which is completely
different from the one presented above but rather extends the 
algebraic geometry approach from
\cite{AW},
have found a nice set of inequalities defining the polytope
$\Delta^G_l$
(where, as before,
$G$
is a compact complex semi-simple connected and 
simply-connected Lie group).
Their inequalities are written in intrinsic 
terms of the structure of the Lie algebra of
$G$.

\bigskip
\begin{rem}
\label{rem-general-abw-ineq-any-homom-with-finite-kernel}
{\rm
Observe that if
$G$
is a finite-dimensional connected
Lie 
\break
group then one has an analogue of 
Theorem~\ref{thm-abw-ineqs}
for any homomorphism   
$G\to {\hbox{\it Ham}}\, (M,\omega)$
with a finite kernel. 
As before this circumstance allows us to pull back the Finsler norm on
the tangent bundle of
${\hbox{\it Ham}}\, (M,\omega)$
to the tangent bundle of
$G$,
where it determines a genuine Finsler metric on the group 
$G$
itself.
}
\end{rem}
\smallskip

\vfil
\eject
\bigskip
\section{Fibrations over a surface with boundary, K-area and
weak coupling}
\label{sect-fibrs-k-area-weak-coupling}

\bigskip
\subsection{Our favorite surface
$\Sigma$ 
}
\label{subsect-surface}

Let
$\Sigma$
be a compact connected oriented Riemann surface of genus 
$0$
with
$l\geq 1$
boundary components:
$\partial\Sigma = T_1\sqcup\ldots\sqcup T_l$.
Fix a volume form 
$\Omega$
on
$\Sigma$
so that
$\int_\Sigma \Omega = 1$.
According to the Moser's theorem
\cite{Mos}, 
any two such volume forms coinciding near
$\partial\Sigma$
can be mapped into each other by a diffeomorphism of
$\Sigma$. 

The orientation of
$\Sigma$
determines an orientation on each
$T_i$
(rotate an outward normal vector to
$\Sigma$
by 90 degrees counterclockwise to get the positive direction on
$\partial\Sigma$).

\bigskip
\subsection{Fibrations, connections, curvatures and holonomies}
\label{subsect-fibr-conn-curv-def}

Let 
$G$
be a connected Lie group
whose tangent bundle is equipped with a bi-invariant Finsler norm
that defines a pseudo-metric
$\rho$
on
$G$. 
Identify the Lie algebra 
${\mathfrak g}$
of
$G$
with the space of right-invariant vector fields on
$G$.
Below we will not use any results from the Lie theory and thus all our
considerations will hold even if
$G$
is an infinite-dimensional Lie group.
Suppose that
$G$
acts effectively on a connected manifold
$F$
(i.e. one has a monomorphism
$G\to {\hbox{\it Diff}}\, (F)$).
The case which is most important for us is when
$F = (M,\omega)$
is a symplectic manifold and
$G= {\hbox{\it Ham}}\, (M,\omega)$. 

Consider the trivial 
$G{\hbox{\rm -bundle}}$
$\pi: P\to\Sigma$, 
with the fiber
$F$.
Let us consider
$G{\hbox{\rm -connections}}$
on the bundle
$P\to\Sigma$,
i.e. the connections whose parallel transports belong to the structural group
$G$.
Let
$L^\nabla$
denote the curvature of a connection 
$\nabla$
on the bundle
$\pi: P\to\Sigma$.
If the fiber
$\pi^{-1} (x)$
is identified with
$F$
then to a pair of vectors
$v,w\in T_x \Sigma$
the curvature tensor associates an element
$L^\nabla (v,w)\in {\mathfrak g}$.
Here we use the fact that
$G$
acts effectively on
$F$.
If no identification of
$\pi^{-1} (x)$
with
$F$
is fixed then
$L^\nabla (v,w)\in {\mathfrak g}$
is defined up to the adjoint action of
$G$
on
${\mathfrak g}$.
Thus if
$\|\cdot \|$
is the bi-invariant Finsler norm on
${\mathfrak g}$
defining our Finsler norm
on
$T_\ast G$
then
$\| L^\nabla (v,w) \|$
does not depend on the identification of
$\pi^{-1} (x)$
with
$F$.

\bigskip 
\begin{defin}
\label{def-norms}
{\rm
We define
$\|L^\nabla\|$
as
\[
{\|L^\nabla\|} = \max_{v,w}
\frac{\| L^\nabla (v,w)\| }{|\Omega (v,w) |},
\]
where the maximum is taken over all pairs
$(v,w)\in T_\ast \Sigma \times T_\ast \Sigma$
such that
$\Omega (v,w)\neq 0$.
}
\end{defin}
\smallskip

\bigskip
\subsection{The definition of K-area}
\label{subsect-def-k-area}

Given a 
$G{\hbox{\rm -connection}}$
$\nabla$
on
$\pi: P\to\Sigma$
the holonomy of  
$\nabla$
along a loop based at
$x\in\Sigma$
can be viewed as an element of
$G$
acting on 
$F$
provided that the bundle is trivialized over 
$x$.
If the trivialization of the bundle is allowed to vary 
then the holonomy is defined up to conjugation in the group
$G$.
Observe that the action of the gauge group does not change
$\|L^\nabla\|$.

Now let
${\cal C}=({\cal C}_1,\ldots,{\cal C}_l)$
be some conjugacy classes in 
$G$.

\bigskip 
\begin{defin}
\label{def-cal-L}
{\rm
Let  
${\cal L} ({\cal C})$
denote the set of connections
$\nabla$
on
$P\to \Sigma$
which are flat over a neighborhood of
$\partial\Sigma$
and such that for any
$i=1,\ldots,l$
and the holonomy of 
$\nabla$
along the oriented boundary component
$T_i\subset \partial\Sigma$
lies in
${\cal C}_i$.
}
\end{defin}
\smallskip

\bigskip 
\begin{defin}
\label{def-K-area}
{\rm
The number 
$0< {\hbox{\it K-area}}\, ({\cal C})\leq +\infty$
is defined as 

\begin{equation}
\label{eqn-def-k-area}
{\hbox{\it K-area}}\, ({\cal C}) =
\sup_{\nabla\in {\cal L} ({\cal C}) }\|L^\nabla\|^{-1}.
\end{equation}

}
\end{defin}
\smallskip

Obviously, since none of the boundary components of
$\Sigma$
is preferred over the others, the quantity
${\hbox{\it K-area}}\, ({\cal C}_1,\ldots, {\cal C}_l)$
does not depend on the order of the conjugacy classes
${\cal C}_1,\ldots, {\cal C}_l$.

\bigskip
\subsection{The relation between K-area and
$\Upsilon_l$}
\label{subsect-thm-k-area-upsilon}

The principal relation between
K-area and
$\Upsilon_l$
is expressed in the following theorem.

\bigskip
\begin{thm}
\label{thm-k-area-distance}
\[
\Upsilon_l\, ({\cal C}) 
=
\frac{1}{{\hbox{\it K-area}}\, ({\cal C})}.
\]
(If K-area is infinite its inverse is assumed to be zero).
\end{thm}
\smallskip

We will prove
Theorem~\ref{thm-k-area-distance}
in a stronger form
(see 
Theorem~\ref{thm-k-area-distance-hom-classes})
but in order to formulate it we need
to introduce more definitions.

\bigskip
\subsection{Systems of paths and their homotopy classes}
\label{subsect-systems-of-paths}

The quantity 
$\Upsilon_l$
can be defined in a different way as follows.

\bigskip 
\begin{defin}
\label{def-system-of-paths}
\rm{
A
{\it
system of paths 
}
$a = (a_1,\ldots, a_l)$
is a tuple of some smooth paths 
$a_1,\ldots, a_l: [0,1]\to G$
such that 
\[
a_1 (0)\cdot\ldots\cdot a_l (0) = Id.
\]
}
\end{defin}
\smallskip

The 
${\hbox{\it length}}\, (a)$ 
of a system of paths
$a$
is defined as the sum of lengths of the paths that form the system, where
the length of a path is measured with respect to the fixed Finsler 
pseudo-metric
on
$G$.

\bigskip 
\begin{defin}
\label{def-system-of-paths-conj-classes}
\rm{
Let
${\cal C} = ({\cal C}_1,\ldots,{\cal C}_l)$
be some conjugacy classes in
$G$.
We define 
${\cal G} ({\cal C})$
as the set of all systems of paths 
$(a_1,\ldots, a_l)$
such that
$a_i (1)\in {\cal C}_i$,
$i=1,\ldots,l$.
}
\end{defin}

\bigskip 
\begin{prop}
\label{prop-upsilon-via-systems-of-paths}
\begin{equation}
\label{eqn-upsilon-syst-paths}
\Upsilon_l\, ({\cal C}) = 
\inf_{\scriptscriptstyle a\in {\cal G} \, ({\cal C})}
{\hbox {\it length}}\, (a) .
\end{equation}
\end{prop}
\smallskip

We will prove
Proposition~\ref{prop-upsilon-via-systems-of-paths}
in 
Section~\ref{sect-pf-prop-upsilon-via-systems-of-paths}.

In the case when the group
$G$
is not simply connected
the space 
${\cal G} ({\cal C})$ 
might have more than one connected component: 
there might be systems of paths in
${\cal G} ({\cal C})$ 
that are not homotopic to each other.
Denote such a homotopy class of a system of paths
$a$
by
$[a]$.
and denote the corresponding connected component of 
${\cal G} ({\cal C})$ 
by
${\cal G}_{[a]}\, ({\cal C})$.

\bigskip
\begin{defin}
\label{def-upsilon-hom-classes-via-system-of-paths}
{\rm
For a homotopy class 
$[a]$
of systems of paths from
${\cal G} ({\cal C})$ 
define
$\Upsilon_{l,[a]}\, ({\cal C})$
by taking in
(\ref{eqn-upsilon-syst-paths})
the infimum only over the systems of paths from
${\cal G}_{[a]}\, ({\cal C})$.
}
\end{defin}
\smallskip
Thus
$\Upsilon_l ({\cal C}) = \inf_{[a]} \Upsilon_{l, [a]} ({\cal C})$.

Just as the space
${\cal G} ({\cal C})$ 
of systems of paths may not be connected,
the space of connections 
${\cal L} ({\cal C})$
on the trivial principal
$G{\hbox{\rm -bundle}}$
$P\to\Sigma$
also might have many connected components, i.e. different connections
from
${\cal L} ({\cal C})$
might not be homotopic to each other in
${\cal L} ({\cal C})$.
Such a homotopy class of a connection 
$\nabla$
will be denoted by
$[\nabla]$
and the corresponding connected component of
${\cal L} ({\cal C})$
will be denoted as
${\cal L}_{[\nabla ]} ({\cal C})$.
If one fixes a homotopy class
$[\nabla ]$
one can define
${\hbox{\it K-area}}_{[\nabla]}\, ({\cal C})$
by taking  the supremum in 
(\ref{eqn-def-k-area}) 
over
${\cal L}_{[\nabla ]} ({\cal C})$.

\bigskip
\subsection{A stronger version of
Theorem~\ref{thm-k-area-distance}}
\label{subsect-str-version-of-thm-k-area-distance}

Now we are able to state a stronger version of
Theorem~\ref{thm-k-area-distance}.

\bigskip
\begin{thm}
\label{thm-k-area-distance-hom-classes}
To a homotopy class 
$[\nabla]$
of connections from
${\cal L} ({\cal C})$
one can naturally associate in a surjective way a homotopy class
$[[\nabla ]] = [a] ([\nabla ])$
of systems of paths from
${\cal G} ({\cal C})$ 
so that
\[
\Upsilon_{l, [[\nabla]]}\, ({\cal C}) 
= 
\frac{1}{{\hbox{\it K-area}}_{[\nabla]}\, ({\cal C})}.
\]
\end{thm}
\smallskip

Theorem~\ref{thm-k-area-distance-hom-classes}
will be proven in 
Section~\ref{sect-proof-thm-k-area-distance}.
As it will easily follow from the proof the correspondence
$[\nabla ]\to [[\nabla ]]$
between homotopy classes of connections and homotopy classes of systems of paths
satisfies the relation below.

For each of the conjugacy classes
${\cal C}_i$,
$i=1,\ldots ,l$,
consider all paths in our group 
$G$
connecting 
${\cal C}_i$
with the identity.
Pick a homotopy class
$\{ c_i\}$
of such paths.
The homotopy class
$\{ c_i\}$
determines a certain conjugacy class
$\tilde{\cal C}_i$
in the universal cover
$\tilde{G}$
of
$G$.

Given the homotopy classes
$\{ c\} = (\{ c_1\} ,\ldots ,\{ c_l\} )$
as above and a system of paths 
$a= (a_1, \ldots , a_l)$
from
${\cal G} ({\cal C})$
complete each path
$a_i: [0,1]\to G$
with
$a_i (1)\in {\cal C}_i$
by a curve connecting
$a_i (1)$
with the identity and representing a path from 
$\{ c_i\}$.
In this way one gets a curve
$\tilde{a}_i: [0,1]\to G$
such that
$\tilde{a}_i (1) = Id$,
$\tilde{a}_i (0) = a_i (0)$.
Since, according to the definition of systems of paths,
$a_1 (0)\cdot\ldots\cdot a_l (0) = Id$,
the pointwise group product 
$t\to \tilde{a}_1 (t)\cdot\ldots\cdot \tilde{a}_l (t)$,
$0\leq t\leq 1$,
of paths
$\tilde{a}_1,\ldots ,\tilde{a}_l$
represents a loop in
$G$
based at the identity. The homotopy class of this loop depends only on
$\{ c\}$
and on the homotopy class
$[a]$
of the system of paths.
If the loop is contractible we say that
$[a]$
{\it
fits}
with
$\{ c\}$.

Denote by
${\cal L}_{\{c\} }^{fit} ({\cal C})$
the set of all the connections 
$\nabla$
in
${\cal L} ({\cal C})$
such that
$[[\nabla ]]$
fits with
$\{ c\}$.
The property of lying in
${\cal L}_{\{c\} }^{fit} ({\cal C})$
depends only on the homotopy class
$[\nabla ]$
of a connection
$\nabla$.

\bigskip
\begin{prop}
\label{prop-hom-classes-conn-id-with-conj-classes-syst-paths-conns}
Let
$\{ c\}$,
${\cal C}$
and
$\tilde{{\cal C}}$
be as above. 
Consider the pseudo-metric on
$\tilde{G}$
defined by the pullback of the Finsler norm from
$T_\ast G$
and let 
$\tilde{\Upsilon}_l\, (\tilde{{\cal C}} )$
be defined by means of this pseudo-metric. 
Then
\[
\tilde{\Upsilon}_l\, (\tilde{{\cal C}} ) = 
\inf_{[[\nabla]]} \, \Upsilon_{l,[[\nabla ]]}\, ({\cal C} ) 
= \inf_{[\nabla ]} \,
\frac{1}{{\hbox{\it K-area}}_{[\nabla]}\, ({\cal C})},
\]
where the infimums are taken over all homotopy classes
$[\nabla ]$
of connections lying in
${\cal L}_{\{c\} }^{fit} ({\cal C})$.

\end{prop}
\smallskip

In particular, this means that in order to study
$\tilde{\Upsilon}_l (\tilde{{\cal C}} )$
by means of K-area we only need to consider connections from
${\cal L}_{\{c\} }^{fit}$.

\bigskip
\subsection{Connections on Hamiltonian fibrations and weak coupling}
\label{subsect-ham-fibr-weak-coupling}

Consider now the specific case of the trivial 
$G{\hbox{\rm -bundle}}$
$P\to \Sigma$
with the fiber
$F = (M,\omega )$,
where
$G= {\hbox{\it Ham}}\,  (M,\omega)$
and
$\Sigma$
is the compact surface of genus 0 with 
$l$
boundary components as above.

Let us briefly recall the following basic definitions (see
\cite{GLS}
for details).

\bigskip 
\begin{defin}
\label{def-Ham-fibr}
{\rm
A closed 2-form
$\tilde{\omega}$
on the total space
$P$
of the bundle
$P\to \Sigma$
is called
{\it 
fiber compatible
}
if its restriction on each fiber of
$P\to \Sigma$
is
$\omega$.
}
\end{defin}
\smallskip

Let us trivialize the bundle
$P=\Sigma\times M\to M$
and let
$pr_M :P\to M$
be the natural projection.
The 
{\it
weak coupling construction}
\cite{GLS}
prescribes that for any
fiber compatible
form
$\tilde{\omega}$
that coincides with
$pr_M^\ast \omega$
near
${\partial P}$
and for sufficiently small
$\varepsilon > 0$
there exists a smooth family of closed 2-forms
$\{ \Omega_\tau \}$,
$\tau \in [0 , \varepsilon ]$,
on
$P$
with the following properties:

\medskip
\noindent
(i) 
$\Omega_0 = \pi^\ast \Omega$,
where
$\pi$
is the projection
$\pi: P\to\Sigma$ 
and
$\Omega$
is the fixed symplectic form on the surface
$\Sigma$;

\smallskip
\noindent
(ii) 
$[\Omega_\tau ] = \tau [\tilde{\omega}] + [\pi^\ast \Omega]$,
where the cohomology classes are taken in
\break
$H^2 (P,\partial P)$;

\smallskip
\noindent
(iii) the restriction of
$\Omega_\tau$
on each fiber of
$\pi$
is a multiple of the symplectic form on that fiber;

\smallskip
\noindent
(iv)
$\Omega_\tau$
is symplectic for
$\tau\in (0,\varepsilon ]$.
\bigskip

\bigskip 
\begin{defin}
\label{def-size}
{\rm
We define
${\hbox{\it size}}\, (\tilde{\omega})$
as the supremum of all
$\varepsilon$
that admit a family
$\{ \Omega_\tau \}$,
$\tau\in [0 , \varepsilon]$,
satisfying the properties (i)-(iv) listed above.
}
\end{defin}
\smallskip

Any 
fiber compatible
form
$\tilde{\omega}$
defines a connection 
$\nabla$
on
$\pi: P\to \Sigma$
and, conversely, any Hamiltonian connection on
$P\to\Sigma$
can be defined by a unique fiber compatible
2-form
$\tilde{\omega}_\nabla$
such that the 2-form on
$\Sigma$
obtained from
$\tilde{\omega}_\nabla^{n+1}$
by fiber integration is 0.
(see e.g.
\cite{GLS}).

Recall that the Lie algebra of
${\hbox{\it Ham}}\,  (M,\omega)$,
which is the algebra of all (globally) Hamiltonian vector fields on
$(M,\omega)$,
is identified with the space of functions on
$(M,\omega)$
with the zero mean value
(see 
Section~\ref{subsect-intro-main-result}).
Therefore the curvature of
$\nabla$
can be viewed as a 2-form associating to each pair 
$v,w\in T_x\Sigma$
of tangent vectors 
on the base a normalized Hamiltonian function 
$H_{v,w}$
on the fiber
$\pi^{-1} (x)$.
The form
$\tilde{\omega}_\nabla$
restricted on the horizontal lifts of vectors 
$v,w\in T_x\Sigma$
at a point
$y\in \pi^{-1} (x)$ 
coincides with 
$H_{v,w} (y)$
(see e.g. 
\cite{GLS}).

Let
$H= (H_1,\ldots , H_l)$
be some (time-dependent) normalized Hamiltonians on
$M$.
Let
$[\varphi_H ] = ([\varphi_{H_1}], \ldots , [\varphi_{H_l}])$
be the conjugacy classes in 
${\hbox{\it Ham}}\,  (M,\omega)$
containing the Hamiltonian symplectomorphisms
$\varphi_{H_1}, \ldots , \varphi_{H_l}$
generated by
$H_1,\ldots , H_l$.
Also let
${\cal C}_H = ({\cal C}_{H_1}, \ldots , {\cal C}_{H_l} )$ 
be the conjugacy classes in
${\widetilde{\hbox{\it Ham}}}\, (M,\omega)$
corresponding to
$H= (H_1,\ldots , H_l)$.

Define 
${\cal F} ([\varphi_H])$
as the set of all the forms 
$\tilde{\omega}_\nabla$,
$\nabla\in {\cal L} ([\varphi_H])$.
Given a homotopy class
$[\nabla ]$
of connections from
$\nabla\in {\cal L} ([\varphi_H])$
define
${\cal F}_{[\nabla ]} ([\varphi_H])$
as the set of all the forms 
$\tilde{\omega}_\nabla$,
$\nabla\in {\cal L}_{[\nabla ]} ([\varphi_H])$.
The Hamiltonian flow generated by
$H_i$,
$1\leq i\leq l$,
represents a homotopy class
$\{ c_i\}$
of paths connecting
$Id$
with
$[\varphi_{H_i} ]$.
Let
$\{ c\} = (\{ c_1\}, \ldots ,\{ c_l\} )$.
In
Section~\ref{subsect-str-version-of-thm-k-area-distance}
we defined the set 
${\cal L}_{\{ c\} }^{fit} ([\varphi_H])$.
Define
${\cal F}_H ([\varphi_H])$
as the set of all
the forms 
$\tilde{\omega}_\nabla$,
$\nabla\in {\cal L}_{\{ c\} }^{fit} ([\varphi_H])$.

\bigskip 
\begin{defin}
\label{def-sigma-weak-coupl}
{\rm
Define the numbers
\[
0< {\hbox{\it size}}\, ([\varphi_H ]),
{\hbox{\it size}}_{[\nabla]}\, ([\varphi_H ]),
{\hbox{\it size}}\, (H)
\leq +\infty
\]
as 
\[
{\hbox{\it size}}\, ([\varphi_H ]) = 
\sup_{\tilde{\omega}\in {\cal F} ([\varphi_H])} 
{\hbox{\it size}}\, (\tilde{\omega}),
\]
\[
{\hbox{\it size}}_{[\nabla]}\, ([\varphi_H ]) = 
\sup_{\tilde{\omega}\in {\cal F}_{[\nabla]} ([\varphi_H])} 
{\hbox{\it size}}\, (\tilde{\omega}),
\]
\[
{\hbox{\it size}}\, (H) = 
\sup_{\tilde{\omega}\in {\cal F}_H ([\varphi_H])} 
{\hbox{\it size}}\, (\tilde{\omega}).
\]
}
\end{defin}
\smallskip 

The following theorem can be proved by exactly the same arguments as the
similar theorems in
\cite{Pol1},
\cite{Pol3}.

\begin{thm}
\label{thm-K-area-coupl-Ham}
Let
${\hbox{\it K-area}}$
be measured, as before, with respect to the norm
$\| h \| = \max_M | h |$
on the Lie algebra of
${\hbox{\it Ham}}\,  (M,\omega)$
(both for the group
${\hbox{\it Ham}}\,  (M,\omega)$
and its universal cover).

Then
\[
{\hbox{\it K-area}}_{[\nabla]}\, ([\varphi_H ])\leq 
{\hbox{\it size}}_{[\nabla]}\, ([\varphi_H]).
\]
for any
$[\nabla]$
and, in general,
\[
{\hbox{\it K-area}}\, ([\varphi_H])\leq 
{\hbox{\it size}}\, ([\varphi_H]).
\]
Also
\[
{\hbox{\it K-area}}\, ({\cal C}_H)\leq 
{\hbox{\it size}} (H).
\]

\end{thm}
\smallskip

\vfil
\eject

\bigskip
\section{The classes
${\cal T}_\tau (H)$,
${\cal T} (H)$,
${\cal T}^0 (H)$
of almost complex structures 
and moduli spaces
${\cal M} 
(\hat{\gamma}, H,\tilde{J})$}
\label{sect-psh-curves-gw-invariants}

In this section we will fix a trivialization
$P = \Sigma\times M$
of trivial bundle 
$P\to\Sigma$
that we considered before.
We rescale our surface
$\Sigma$
to present it as a non-compact surface of area
$1$
with
$l$
cylindrical ends. We fix an identification 
$\Phi_i:  [0,+\infty)\times S^1\to \Sigma_i$,
$1\leq i\leq l$,
of each end
$\Sigma_i\subset \Sigma$
with the standard cylinder
$[0,+\infty) \times S^1$.
Without loss of generality we may assume that the identifications
are chosen in such a way that near infinity 
the conformal structure on the ends gets identified
with the standard conformal structure on the cylinder
$[0,+\infty) \times S^1$.

If
$J$
is an almost complex structure on
$M$
compatible with 
$\omega$
and
$\tilde{J}$
is an almost complex structure
$\tilde{J}$
on 
$\Sigma\times M$
we say that
$\tilde{J}$
is
$J{\hbox{\it -fibered}}$
if the following conditions are fulfilled:

\medskip
\noindent
$\bullet$
$\tilde{J}$
preserves
the tangent spaces to the fibers of
$\pi : \Sigma\times M\to \Sigma$;

\smallskip
\noindent
$\bullet$
the restriction of 
$\tilde{J}$ 
on any fiber is an almost complex structure compatible with the
symplectic form
$\omega$
on the fiber; 

\smallskip    
\noindent     
$\bullet$
the restriction of
$\tilde{J}$ 
on any fiber
$\pi^{-1} (x)$
for
$x$
outside of some compact subset of
$\Sigma$
is 
$J$.
\medskip

Now let
$H = (H_1,\ldots, H_l)$
be some (time-dependent) Hamiltonians on
$M$.
Let us choose some
${\hat{\gamma}}\in {\cal P} (H)$.
Let us also pick once and for all a cut-off function
$\beta: {\bf R}\to [0,1]$
such that
$\beta (s)$
vanishes for
$s\leq\epsilon$
and
$\beta (s) = 1$
for
$s\geq 1-\epsilon$
for some small
$\epsilon>0$.

Any section
$u:\Sigma\to P$
by means of the trivialization 
$P =\Sigma\times M$
induces some maps
\[
u_i = u\circ \Phi_i : [0,+\infty)\times S^1\to M.
\]
Suppose that
$J$
is an
$\omega{\hbox{\rm -compatible}}$
almost complex structure on
$M$. 
For each 
$i = 1,\ldots,l$
consider the non-homogeneous Cauchy-Riemann equation
\begin{equation}
\label{eqn-nonhom-cauchy-riem}
\partial_s u_i + J (u_i)\partial_t u_i - \beta (s) \nabla_u H_i (t,u_i) 
= 0,
\end{equation}
where gradient is taken with respect to the Riemannian metric 
$g (\cdot ,\cdot ) =\omega (\cdot , J\cdot )$
on
$M$.
According to
\cite{Gro-pshc}, 
$1.4.C$,
the solutions of such an equation correspond exactly to the
pseudo-holomorphic sections of
$\pi^{-1} (\Sigma_i)\to \Sigma_i$
with respect to some unique 
$J{\hbox{\rm -fibered}}$ 
almost complex structure on
$\pi^{-1} (\Sigma_i)$
in the following way.
Fix an almost complex structure
$j$
on
$\Sigma$
which is compatible with the symplectic form
$\Omega$
and which extends the standard complex structures on the cylinders
$\Sigma_i = \Phi_i ([0,+\infty)\times S^1 )$,
$i=1,\ldots ,l$.

\bigskip 
\begin{defin}
\label{def-slow-acs-on-a-bundle}
{\rm
Let
$\tilde{J}$
be an almost complex structure on
$P$
and let
$H = (H_1,\ldots , H_l)$
be Hamiltonians as above.
We shall say that
$\tilde{J}$
is
{\it
$H{\hbox{\it -compatible}}$
}
if there exists  an almost complex structure
$J = J (\tilde{J})$
on
$M$
compatible with
$\omega$
such that the following conditions hold.

\smallskip
\noindent
$\bullet$
$\tilde{J}$
is 
$J{\hbox{\rm -fibered}}$.

\smallskip
\noindent
$\bullet$
$\pi\circ \tilde{J} = j\circ\pi$,
where 
$\pi: \Sigma\times M\to\Sigma$
is the projection.

\smallskip
\noindent
$\bullet$
For each
$i=1,\ldots, l$
there exist some number
$K_i$
such that over
${\tilde{\Sigma }}_i = \Phi_i ([K_i, +\infty )\times S^1) \subset
{\Sigma }_i = \Phi_i ([0, +\infty )\times S^1 )$
the
$\tilde{J}{\hbox{\rm -holomorphic}}$
sections of
$\pi^{-1} ({\tilde{\Sigma }}_i)\to {\tilde{\Sigma }}_i$,
viewed as maps
$u_i: [K_i, +\infty )\times S^1 \to M$,
are exactly the solutions of the non-homogeneous Cauchy-Riemann equation
(\ref{eqn-nonhom-cauchy-riem})
for 
$J$.
\smallskip

Let us denote by
${\cal T} (H)$
the space of all 
$H{\hbox{\rm -compatible}}$
almost complex structures on
$\Sigma\times M$.
We also define a subset
${\cal T}^0 (H)\subset {\cal T} (H)$
that includes only those 
$\tilde{J}\in {\cal T} (H)$
which for some
$J$
are split as
$\tilde{J} = j\times J$
over a compact part of
$\Sigma$
that contains
$\Sigma\setminus \bigcup_{i=1}^l \Sigma_i$,
where
$J = J ({\tilde{J}})$.
In such a case we will denote
$\tilde{J} = {\tilde{J}}_{J, H}$.
}
\end{defin}
\medskip

Let
$[\varphi_H] = ([\varphi_{H_1} ],\ldots ,[\varphi_{H_l}])$
be the conjugacy classes in
${\hbox{\it Ham}}\,  (M,\omega)$
as before.

Now we define the classes
${\cal T}_\tau (H)$
of almost complex structures that are used in the definition of a durable
element in
${\cal P} (H)$
(see
Section~\ref{subsect-intro-main-result}).

\bigskip
\noindent
\begin{defin}
\label{def-cal-T-tau}
{\rm
Consider all the families
$\{ \Omega_{\tilde{\omega}_\nabla, \tau} \}$,
that arise from the weak coupling 
construction associated with
$\tilde{\omega}_\nabla$,
$\nabla\in {\cal L}_{\{ c\} }^{fit} ([\varphi_H])$
(see
Section~\ref{subsect-ham-fibr-weak-coupling}).
Given a number
$\tau_0\in (0,  {\hbox{\it size}}\, (H) )$
consider the set
${\cal Q}_{\tau_0}$ 
of all the forms
$\Omega_{\tilde{\omega}_\nabla, \tau_0}$
from the families
$\{ \Omega_{\tilde{\omega}_\nabla, \tau} \}$
as above
(i.e. we consider only those families which are defined for
the value
$\tau_0$
of the parameter
$\tau$
and pick the form
$\Omega_{\tilde{\omega}_\nabla, \tau_0}$
from each such family).
Denote by
${\cal T}_{\tau_0} (H)$
the set of all the almost complex structures in
${\cal T} (H)$
which are compatible with some symplectic from
${\cal Q}_{\tau_0}$.
}
\end{defin}
\smallskip

Below we also give a precise definition of the moduli spaces
${\cal M} 
(\hat{\gamma}, H, \tilde{J})$
involved in the definition of a durable
element in
${\cal P} (H)$
(see
Section~\ref{subsect-intro-main-result}).

Let
$\hat{\gamma} = [{\hat{\gamma}}_1,\ldots ,{\hat{\gamma}}_l ]$,
where
${\hat{\gamma}}_i = [\gamma_i, f_i ]\in {\cal P} (H_i)$,
$i=1,\ldots ,l$.
Let
$\tilde{J}\in {\cal T} (H)$.

\bigskip 
\begin{defin}[cf. \cite{PSS}]
\label{def-modul-space}
{\rm
Denote by
$
{\cal M} 
(\hat{\gamma}, H, \tilde{J})
$
the space of all smooth 
$\tilde{J}{\hbox{\rm -holomorphic}}$
sections
$u:\Sigma \to P$
which satisfy the following conditions.

\medskip
\noindent
(i) The maps
$u_i = u\circ \Phi_i : [0,+\infty)\times S^1\to M$
constructed above satisfy 

\[ 
\gamma_i (t) = \lim_{s\to +\infty} u_i (s,t).
\]

\smallskip
\noindent
(ii) The closed surface obtained by capping off  
$u(\Sigma)\subset M$
with the discs
$f_1,\ldots ,f_l$
(taken with the opposite orientations)
represents a torsion homology class in
$H_2 (M,{\bf Z})$.

}
\end{defin}
\smallskip

Looking at the condition (ii) one sees that the space 
${\cal M} 
(\hat{\gamma}, H,\tilde{J})$
depends only on the equivalence class
$\hat{\gamma} = [{\hat{\gamma}}_1,\ldots ,{\hat{\gamma}}_l ]$
of the tuple
$({\hat{\gamma }}_1,\ldots,{\hat{\gamma }}_l )$
-- see the definition of the equivalence relation in
Section~\ref{subsect-intro-main-result}.

\vfil
\eject

\bigskip
\bigskip
\section{Proof of 
Theorem~\ref{thm-exist-psh-curve-implies-estimate}
} 
\label{sect-proof-thm-exist-psh-curve-implies-estimate}

Let
$\Upsilon_l^{Ham}$
be defined on
${\hbox{\it Ham}}\,  (M,\omega)$
using the Hofer metric.
Let
$\Sigma$
be the surface with
$l$
cylindrical ends, as before.
Let
$\{ c\} = (\{ c_1\}, \ldots ,\{ c_l\} )$
be as in
Section~\ref{sect-psh-curves-gw-invariants}.

From 
Proposition~\ref{prop-hom-classes-conn-id-with-conj-classes-syst-paths-conns}
one gets that
\begin{equation}
\label{eqn-upsilon-ham-covering}
\Upsilon_l\, ({\cal C}_H) = 
\inf_{\nabla} \Upsilon_{l,[[\nabla ]]}^{Ham}\, ([\varphi_H])
\end{equation}
where the infimum is taken over all the connections
$\nabla\in {\cal L}_{\{ c\} }^{fit} ([\varphi_H])$.

It follows from the hypothesis of the theorem that for any
$\tau_0$,
$0<\tau_0 < {\hbox{\it size}}\, (H)$, 
there exist a connection
$\nabla\in {\cal L}_{\{ c\} }^{fit} ([\varphi_H])$,
the corresponding 2-form
$\tilde{\omega}_\nabla$
such that
$\tau_0 < {\hbox{\it size}}\, (\tilde{\omega}_\nabla )$,
and a weak coupling deformation
$\{\Omega_{\tilde{\omega}_\nabla , \tau }\}$,
$0< \tau < {\hbox{\it size}}\, (\tilde{\omega}_\nabla )$,
so that the space
${\cal M} (\hat{\gamma}, H, \tilde{J}_{\tau_0} )$
is non-empty for some
$\tilde{J}_{\tau_0} \in {\cal T}_{\tau_0} (H)$
compatible with the symplectic form
$\Omega_{\tilde{\omega}_\nabla , \tau_0 }$.

Pick a map
$u\in {\cal M}
(\hat{\gamma}, H, \tilde{J}_{\tau_0} )$.
Since
$\tilde{J}_{\tau_0}$
is compatible with the symplectic form
$\Omega_{\tilde{\omega}_\nabla,\tau_0 }$,
one has 
\begin{equation}
\label{eqn-ineq-tot-integral-positive}
0\leq \int_{u(\Sigma)} \Omega_{\tilde{\omega}_\nabla ,\tau_0 } = 
\int_{u(\Sigma)} \tau_0 \tilde{\omega} + 
\int_{u(\Sigma)}\pi^\ast \Omega,
\end{equation}
where
$u(\Sigma)$
is viewed as a surface in
$P$
and
$\pi: P\to\Sigma$
is the projection.

\bigskip
\begin{lem}
\label{lem-intergal-equal-action}
\begin{equation}
\label{equal-action-omega-tilde}
\int_{u(\Sigma)} \tilde{\omega}_\nabla = - {\cal A}_H (\hat{\gamma}).
\end{equation}
\end{lem}
\smallskip

Postponing the proof of the lemma we first finish the proof of the theorem.
Indeed, since the total 
$\Omega{\hbox{\rm -area}}$ of
$\Sigma$
is 1 one can rewrite 
(\ref{eqn-ineq-tot-integral-positive}),
in the case when
${\cal A}_H (\hat{\gamma}) = - \int_{u(\Sigma)}\tilde{\omega}_\nabla$
is positive
(otherwise the theorem is trivial)
as
\[
\tau_0\geq \frac{1}{{\cal A}_H (\hat{\gamma}) }
\]
Since
$\nabla\in {\cal L}_{\{ c\} }^{fit} ([\varphi_H ])$
was chosen arbitrarily,
$\tau_0$
can be chosen arbitrarily close to
${\hbox{\it size}}\, (H)$
one gets that
\[
{\hbox{\it size}}\, (H)
\geq \frac{1}{{\cal A}_H (\hat{\gamma}) }.
\]
Using this inequality together with
(\ref{eqn-upsilon-ham-covering}),
Theorem~\ref{thm-k-area-distance-hom-classes}
and
Theorem~\ref{thm-K-area-coupl-Ham}
one readily obtains the needed result.
\b

%\vfil
%\eject
\bigskip
\bigskip
\noindent
{\bf The proof of
Lemma~\ref{lem-intergal-equal-action}.}

This is a purely topological fact  -- we have already used all the complex
properties of
$u$
that we needed.
Therefore we can rescale 
$\Sigma$
back and make it a compact surface with boundary.
Then one can extend
${\Phi}_i$,
$i=1,\ldots ,l$,
to a map
${\Phi}_i : [0,+\infty]\times S^1\to \Sigma$ 
so that
${\Phi}_i (+\infty \times t)$,
$o\leq t\leq 1$,
parameterizes the boundary component
$T_i$
of
$\Sigma$.
The map 
$u$
restricted on the boundary component
$T_i$
of
$\Sigma$
produces the curve
$\gamma_i\in M$,
$i=1,\ldots,l$.
We cap off the boundaries of
$\Sigma$
with some discs 
$D_1,\ldots, D_l$
and get a closed surface
$\hat{\Sigma}$.

Given a connection
$\nabla\in {\cal L}_{\{ c\} }^{fit} ([\varphi_H])$
on the bundle over the compact surface
$\Sigma$
one can assume without loss of generality that
the trivialization of 
$P\to\Sigma$
is already adjusted in such a way that for any
$\tau$
the holonomy of
$\nabla$
(taken with respect to the trivialization of 
$P\to\Sigma$) 
over the path
${\{ t\to {\Phi}_i (+\infty \times t)\}}_{0\leq t\leq \tau}$
is the flow
$\varphi^\tau_{H_i}$,
$i=1,\ldots, l$,
of
$H_i$
for the time 
$\tau$.

Then it is not hard to prove the following technical 
sublemma (see e.g.
\cite{Sch-1},
Section 4.1).

\bigskip
\begin{sublem}
\label{lem-cap-off}
Over a disc
$D_i$,
$i=1,\ldots,l$,
one can construct a trivialized bundle 
$E_i = D_i\times M\to D_i$
together with its section 
$U_i: D_i\to E_i$
and a connection
$\nabla_i$
with the 2-form
$\tilde{\omega}_{\nabla_i}$
on
$E_i$
with the following properties.

\medskip
\noindent
$(i)$
The trivialized bundle
$E_i\to D_i$
agrees along
$\partial D_i = \Delta_i$
with the trivialized bundle
$P\to\Sigma$.

\smallskip
\noindent
$(ii)$
If by means of the trivialization the section
$U_i: D_i\to E_i$
is viewed as map 
$U_i : D_i\to M$
then
$U_i (D_i) = f_i$,
$U_i (\partial D_i) = \gamma_i$.

\smallskip
\noindent
$(iii)$
The connection
$\nabla_i$
$\tilde{\omega}_i$
on
$E_i$
smoothly extends the connection
$\nabla$
on 
$P\to\Sigma$.

\smallskip
\noindent
$(iv)$
$\int_{D_i} U_i^\ast \tilde{\omega}_{\nabla_i} = 
{\cal A}_{H_i} ([\gamma_i, f_i])$.
\end{sublem}
\smallskip

Note that the definition of action functional in
\cite{Sch-1}
differs from the one we use here by the sign at the term containing the 
integral over a disc
and the Hamiltonians used in
\cite{Sch-1}
actually correspond not to our Hamiltonians
$H = (H_1,\ldots, H_l )$
but to
${\bar H} = ({\bar H}_1,\ldots, {\bar H}_l )$,
due to a different sign convention given by
(\ref{eqn-def-ham-vect-field})
so that at the end we get the plus sign in the right-hand side of
$(iv)$
above.

The lemma allows us the bundles
$E_i\to D_i$, 
the maps
$U_i$, 
and the forms
$\tilde{\omega}_{\nabla_i}$,
$i=1,\ldots,l$,
to 
$P\to\Sigma$, 
$u$
and 
$\tilde{\omega}_\nabla$
respectively. Namely, consider again
$\Sigma$
as a compact surface with boundary and consider a closed surface
$\hat{\Sigma }$
obtained by capping off each boundary component of
$\Sigma$
with a disc. 
We construct a trivialized bundle 
$\hat{P} = \hat{\Sigma}\times M\to\hat{\Sigma}$
to which we extend the section
$u$
(viewed now as a map
$u: \hat{\Sigma}\to \hat{P}$)
and the 2-form
$\tilde{\omega}_\nabla$
(viewed now as a form on
$\hat{P}$).
Recalling the condition 
(ii)
of 
Definition~\ref{def-modul-space}, 
we get
\[
0 = \int_{u(\hat{\Sigma})} \tilde{\omega}_\nabla 
\]
and hence, according to 
Lemma~\ref{lem-cap-off},
\[
\int_{u(\Sigma)} \tilde{\omega}_\nabla =
- \sum_{i=1}^l \int_{U_i (D_i)} \tilde{\omega}_\nabla = 
- {\cal A}_H (\hat{\gamma}).
\]
This gives us the 
equality
(\ref{equal-action-omega-tilde})
and finishes the proof of the lemma.
\b
\smallskip

\vfil
\eject

\bigskip
\bigskip
\section{The proof of 
Theorem~\ref{thm-main-time-indep-hams}}
\label{sect-pf-of-main-thm-slow-hams}

Suppose, as in the hypothesis of 
Theorem \ref{thm-main-time-indep-hams},
that for some singular cohomology classes 
$c_1,\ldots, c_{l-1}\in H^\ast (M, {\bf Q})$
one has
\begin{equation}
\label{eqn-1-identity-to-be-involved-in}
c_1\ast\ldots\ast c_{l-1} = \sum_{B\in\Pi} c_B e^{2\pi i B},
\end{equation}
where
$c_B\in H^\ast (M, {\bf Q})$.

Theorem~\ref{thm-main-time-indep-hams}
follows immediately from
Theorem~\ref{thm-exist-psh-curve-implies-estimate}
and the following proposition proved below.

\bigskip
\begin{prop}
\label{prop-quantum-multipl-implies-nontriv-gw-number}
In the notation as in
Definition~\ref{def-involved}
suppose that
$\hat{\gamma}\in {\cal P} (H)$
associated with 
$z^1,\ldots, z^l, A$
is involved in an identity
(\ref{eqn-1-identity-to-be-involved-in}).
Then 
$\hat{\gamma}\in {\cal P} (H)$
is durable.
\end{prop}
\smallskip

To check that 
$\hat{\gamma}\in {\cal P} (H)$
is durable one needs to study of the multiplicative structure 
of the 
{\it
Floer cohomology 
}
of
$(M,\omega)$
defined below.

\bigskip
\subsection{Conley-Zehnder index of a periodic trajectory of a 
Ham\-il\-tonian flow}
\label{subsect-Conley-Zehnder-index}

Let
$h$
be a Hamiltonian on
$M$.
For an element
$\hat{\gamma}= [\gamma, f]\in {\cal P} (h)$,
one can define its
{\it
Conley-Zehnder index
}
$\mu (\hat{\gamma})\in {\bf Z}$
(see 
\cite{Co-Ze}).
The Conley-Zehnder index satisfies the property
\[
\mu (A\sharp \hat{\gamma}) = \mu (\hat{\gamma}) - 2 c_1 (A).
\]

If 
$h$
is a slow Morse Hamiltonian and 
$\hat{\gamma} = [\gamma, f]$
is formed by constant maps into a critical point 
$y$
of
$H$,
then
$\mu (\hat{\gamma})$
is equal to the Morse coindex of
$y$
(i.e. 
$2n$
minus the Morse index).

For
$H = (H_1,\ldots ,H_l)$
on
$M$
and an element
$\hat{\gamma} = [{\hat{\gamma}}_1, \ldots ,{\hat{\gamma}}_l ]
\in {\cal P} (H)$
we define
\[
\mu (\hat{\gamma}) = \sum_{i=1}^l \mu({\hat{\gamma}}_i).
\]

\bigskip
\subsection{Floer cohomology}
\label{subsect-floer-quantum-basic-facts}

Now we recall some basic facts about Floer cohomology.

Assume that 
$(M,\omega)$,
is a closed connected strongly semi-positive symplectic manifold.
The Floer theory shows that for a Hamiltonian 
$h: S^1\times M\to {\bf R}$
one can define a kind of ``Morse homology''
for the action functional
${\cal A}_h$ 
over the ring
$\Lambda_\omega$.
The construction goes along the following lines
(for the details of the theory in the semi-positive case see
\cite{Ho-Sa}). 
Given a generic (a so called
{\it 
regular
}) 
pair
$(h,J)$,
where
$h$
is a (time-dependent) Hamiltonian and
$J$
is an almost complex structure compatible with
$\omega$,
one defines a chain complex 
$CF_\ast (h,J)$
of
$\Lambda_\omega{\hbox{\rm -modules}}$.
Here 
$CF_\ast (h,J)$
is a free graded module over the graded ring
$\Lambda_\omega$
generated by the critical points of
${\cal A}_h$
which are time-1 contractible periodic orbits of the
Hamiltonian flow of
$h$.
The grading of such an orbit is given by the Conley-Zehnder index.
The differential 
$\partial$
in the complex is defined similarly to the 
finite-dimensional Morse homology by means of 
counting solutions of an appropriate Cauchy-Riemann equation
that represent ``the gradient trajectories''
connecting critical points
of
${\cal A}_h$
of neighboring indices
(see e.g.
\cite{Ho-Sa}
for details).
One can show that
$\partial^2 = 0$
(see
\cite{Ho-Sa}).
The homology group
$HF_\ast (h,J)$
of the chain complex
$CF_\ast (h,J)$
is called
{\it
the Floer homology group
}
(it is actually a module over
$\Lambda_\omega$).
Given two different (regular) pairs
$(h_\alpha, J_\alpha)$
and
$(h_\beta, J_\beta)$
there exists a natural isomorphism
$I^{\beta\alpha}_\ast: 
HF_\ast (h_\alpha, J_\alpha)\to HF_\ast (h_\beta, J_\beta)$.

Assume now that
$h$ 
is a slow Morse Hamiltonian.
Let
$J$
be an almost complex structure on
$M$
compatible with
$\omega$
so that
$h$
is a Morse-Smale function with respect to the Riemannian metric
$g (\cdot ,\cdot ) =\omega (\cdot , J\cdot )$.
The most crucial property of slow Morse Hamiltonians is that for
any such 
$J$
the pair
$(h, J)$
is regular and the Floer chain complex
$CF_\ast (h,J)$
can be identified with the Morse chain complex 
$C_\ast (-h, g)$
(see
\cite{Sal-Ze}).
The seemingly strange combination of signs is due to
our choice of the signs in the definition of action functional
${\cal A}_h$
so that the downward gradient flow of
${\cal A}_h$
corresponds to the upward gradient flow of
$h$.

The appropriate map
$I^{\beta\alpha}_\ast$
provides a natural grading-preserving isomorphism between
$HF_\ast (h,J)$
and the homology of the chain complex
$C_\ast (-h, g)$
which leads to a natural isomorphism
\begin{equation}
\label{eqn-nat-isom-floer-sing-homol}
{\cal SF}_h: H_\ast (M)\otimes \Lambda_\omega \to HF_\ast (h,J)
\end{equation}
such that for any two regular pairs
$(h_\alpha, J_\alpha)$
and
$(h_\beta, J_\beta)$,
with
$h_\alpha, h_\beta$
being slow Morse Hamiltonians, one has
$I^{\beta\alpha}_\ast \circ {\cal SF}_{h_\alpha } = {\cal SF}_{h_\beta}$.

In a similar fashion, given a Hamiltonian
$h: S^1\times M\to {\bf R}$,
one can define the Floer cochain complex
$CF^\ast (h,J) = Hom (CF_\ast (h,J), \Lambda_0)$
whose cohomology 
$HF^\ast (h, J)$ 
is called 
{\it
the Floer cohomology group.
}
There exists a natural Poincar{\'e} isomorphism
$HF^k (h,J) \cong HF_{2n-k} (\bar{h},J)$,
where
$\bar{h} (t,x) = - h (-t,x)$.
If
$h$
is a slow Morse Hamiltonian one can identify the Morse and the Floer cohomology
in the same fashion as homology.

\bigskip
\subsection{The pair-of-pants product in Floer cohomology}
\label{subsect-floer-quantum-multiplication}

We are going to state the basic facts concerning the pair-of-pants product
in Floer cohomology in our specific situation.

Let
$H = (H_1,\ldots, H_l)$
be our slow Morse Hamiltonians.
Let
$J$
be an almost complex structure on
$(M,\omega )$
such that
$H_1,\ldots , H_l$
are Morse-Smale functions with respect to the metric
$g(\cdot, \cdot ) = \omega (\cdot, J\cdot )$.
Given
$0<\tau <{\hbox{\it size}}\, (H)$
denote by
${\cal T}_J (H)$
the set of all
$\tilde{J}\in {\cal T} (H)$
such that
$J = J (\tilde{J})$.
Set
${\cal T}_{\tau, J} (H) = {\cal T}_\tau (H)\cap {\cal T}_J (H)$.

The following statement is similar to the one from
\cite{PSS}
and
\cite{Sch-PhD}
although our setup is slightly different from the one used there:
here we use almost complex structure from
${\cal T} (H)$
while 
in 
\cite{PSS},
\cite{Sch-PhD}
a smaller class
${\cal T}^0 (H)\subset {\cal T}(H)$
is used.

\bigskip
\begin{prop}
\label{prop-moduli-space-generically-smooth}
Assume that
$(M,\omega)$
is strongly semi-positive.
Then for a gene\-ric
$\tilde{J}_\tau \in {\cal T}_{\tau ,J} (H)$,
and any 
$\hat{\gamma}\in {\cal P} (H)$
with
$\mu (\hat{\gamma}) = 2n$
the space
${\cal M} (\hat{\gamma}, H, \tilde{J}_\tau )$
is either empty or an oriented compact zero-dimensional manifold.
\end{prop}
\smallskip

Given such a generic
$\tilde{J}_\tau$
we will say that
the pair
$(H,\tilde{J}_\tau )$
is
{\it
regular.
}

Given a regular pair 
$(H, \tilde{J}_\tau )$,
$\tilde{J}_\tau \in {\cal T}_{\tau ,J} (H)$,
and an element
$\hat{\gamma}\in {\cal P} (H)$
such that
$\mu({\hat{\gamma}}) = 2n$
count the curves from the compact zero-dimensional moduli space
${\cal M} (\hat{\gamma}, H, \tilde{J}_\tau )$
with their signs. 
The resulting Gromov-Witten number
will be denoted by 
$n (\hat{\gamma}, H, \tilde{J}_\tau )$.
Form the sum

\begin{equation}
\label{Floer-cycle}
\theta_{\Sigma, H,\tilde{J}_\tau } = 
\sum_{\hat{\gamma}}
n (\hat{\gamma}, H, \tilde{J}_\tau ) \hat{\gamma},
\end{equation}

\noindent
where the sum is taken over all
$\hat{\gamma}\in {\cal P} (H)$
such that 
$\mu (\hat{\gamma}) = 2n$.

The sum in
(\ref{Floer-cycle})
represents an integral chain in the chain complex
\begin{equation}
\label{eqn-floer-product-chain-complex}
CF_\ast (H_1, J)\otimes\ldots\otimes  CF_\ast (H_l, J).
\end{equation}
Roughly speaking, this chain complex is an integral 
Morse homology complex for the action functional
${\cal A}_H$,
$H = (H_1,\ldots , H_l)$,
and its homology is equal to
$HF_\ast (H_1, J)\otimes\ldots\otimes HF_\ast (H_l, J)$,
where 
$\otimes$
stands for the graded tensor product over the Novikov
ring 
$\Lambda_0$
-- 
see
\cite{PSS}
for details.

The following theorem is a slight generalization of the
main result from
\cite{PSS}:
we use a bigger class of admissible almost complex structures
(as in
Proposition~\ref{prop-moduli-space-generically-smooth}).

\bigskip
\begin{thm}
\label{thm-main-pss}
Assume that
$(M,\omega)$
is strongly semi-positive.
For a regular pair
$(H, \tilde{J}_\tau )$,
$0< \tau < {\hbox{\it size}}\, (H)$,
$\tilde{J}_\tau \in {\cal T}_{\tau ,J} (H)$,
the chain
$\theta_{\Sigma, H,\tilde{J}_\tau }$
defines a cycle in the chain complex
(\ref{eqn-floer-product-chain-complex}). 
The corresponding homology class
\[
\Theta_{H, J}\in 
HF_\ast (H_1, J)\otimes\ldots\otimes HF_\ast (H_l, J)
\]
is of degree
$2n$
and does not depend on
$\tilde{J}_\tau \in {\cal T}_J (H)$.

By means of the Poincar{\'e} duality the homology class
$\Theta_{H,J}$
defines the 
{\rm\bf
pair-of-pants 
product
}
\[
\varrho : 
HF^{i_1} (H_1, J)\otimes\ldots\otimes HF^{i_{l-1}} (H_{l-1}, J) \to 
HF^{i_1+\ldots +i_{l-1}} ({\bar{H}}_l, J)
\] 
on the Floer cohomology 
which is related to the quantum product in the following way.
Given a Hamiltonian function 
$h:S^1\times M\to {\bf R}$
there exists a natural isomorphism
${\cal QF}_h: QH^\ast (M)\to HF^\ast (h,J)$
which intertwines the quantum product on the quantum cohomology with 
the pair-of-pants product 
$\varrho$
on the Floer cohomology: 
\[
{\cal QF}_{{\bar{H}}_l} (b_1\ast\ldots\ast b_{l-1}) = 
\varrho ({\cal QF}_{H_1} (b_1)\otimes\ldots
\otimes {\cal QF}_{H_{l-1}} (b_{l-1})).
\]

\end{thm}
\smallskip

The following statement, I believe, has been known to the experts but 
I was unable
to find it in a published form.

\bigskip
\begin{prop}
\label{prop-isom-betw-Morse-Floer-sing-hom-for-slow-Hams}
If 
$h$
is a slow Morse Hamiltonian then
${\cal QF}_h = {\cal SF}_h$
where
${\cal SF}_h: H^\ast (M)\otimes\Lambda_\omega\to HF^\ast (h,J)$
is the natural isomorphism in cohomology induced by the isomorphism from
(\ref{eqn-nat-isom-floer-sing-homol}).
\end{prop}
\smallskip

\bigskip
\subsection{Proof of
Proposition~\ref{prop-quantum-multipl-implies-nontriv-gw-number}}
\label{subsect-non-vanishing-of-N}

Let the setup be as in
Section~\ref{subsect-floer-quantum-multiplication}.
Using the fact that
$H_1,\ldots , H_l$
are Morse-Smale functions with respect to the metric
$g(\cdot, \cdot ) = \omega (\cdot, J\cdot )$
we identify the Floer chain complex
$CF_\ast (H_i, J)$
and the Morse chain complex
$C_\ast ({\bar H}_i)$,
$i=1,\ldots , l$.
We want to study the properties of the pair-of-pants 
product on the level of cochains in the (dual) Morse complexes for the 
slow Morse Hamiltonians.

Pick an arbitrary
$0< \tau < {\hbox{\it size}}\, (H)$
and a regular pair
$ (H,\tilde{J})$,
$\tilde{J}_\tau \in {\cal T}_{\tau ,J} (H)$. 
In view of the identification of Floer and Morse complexes and
Theorem~\ref{thm-main-pss}
the chain
\begin{equation}
\label{eqn-theta-chain-[M]}
\theta_{\Sigma, H,\tilde{J}_\tau } =
\sum_{\mu (\hat{\gamma}) = 2n}
n (\hat{\gamma}, H, \tilde{J}_\tau )\, \hat{\gamma},
\end{equation}
can be considered as representing a homology class 
$\Theta_{H,J}$
in 
\[
H_\ast (C_\ast ({\bar H}_1))\otimes\ldots\otimes 
H_\ast (C_\ast ({\bar H}_l)).
\]
(Recall that all the tensor products are taken over
$\Lambda_0\subset \Lambda_\omega$).
Theorem~\ref{thm-main-pss}
and
Proposition~\ref{prop-isom-betw-Morse-Floer-sing-hom-for-slow-Hams}
imply that
$\Theta_{H,J}$
defines the
{\it
quantum
}
multiplication on 
$H^\ast (M, \Lambda_\omega )\cong H_\ast (C^\ast  ({\bar H}_i) )$,
$i=1,\ldots ,l$.

According to the hypothesis of 
Proposition~\ref{prop-quantum-multipl-implies-nontriv-gw-number},
$\hat{\gamma}\in {\cal P} (H)$
associated with
$z^1,\ldots, z^l, A$
is involved in the identity
(\ref{eqn-identity-to-be-involved-in})
in the quantum cohomology.
For each
$i=1, \ldots, l-1$
denote by
$\{\alpha_j^i\}$,
$1\leq j\leq {\hbox{\rm dim}}\, H_\ast (M, {\bf Q})$,
the corresponding 
$i{\hbox{\rm -friendly}}$
basis of
$H_\ast (M, {\bf Q})$.
Denote by
$z_j^l$,
$j\in {\cal I}$,
all the critical points of
$H_l$,
where
${\cal I}$
is some finite set of indices.
Recall that the Morse complexes are taken with coefficients in the ring
$\Lambda_\omega$
which contains 
${\bf Q}$
and therefore has
{\it
no torsion}.
Thus one can represent the homology class 
$\Theta_{H,\tilde{J}}$
by a chain as follows:
\begin{equation}
\label{eqn-Theta-as-a-chain-in-Morse-complex}
\Theta_{H, J } =
\sum \alpha_{j_1}^1\otimes\ldots\otimes \alpha_{j_{l-1}}^{l-1}
\otimes [z_{j_l}^l, f_{j_l}^l]. 
\end{equation}
for some 
$f_{j_l}^l$
such that
$[z_{j_l}^l, f_{j_l}^l]\in {\cal P} (H_l)$.
Here each rational homology class
$\alpha_{j_i}^i$,
$i=1,\ldots , l-1$,
is a generator from the 
$i{\hbox{\rm -friendly}}$
basis
$\{ \alpha_j^i\}$
and is viewed as a homology class of 
the Morse chain complex of
${\bar H}_i$
{\it
over
}
${\bf Q}$.
To get the expression
(\ref{eqn-Theta-as-a-chain-in-Morse-complex})
we have incorporated all the coefficients from
$\Lambda_\omega$
in 
$\Theta_{H, J }$
into the last factor of the tensor product and we have expanded a product of
a coefficient from
$\Lambda_\omega$
with an element of the
$\Lambda_\omega{\hbox{\rm -module}}$
$C_\ast ({\bar H}_l)$
as a sum of elements from
${\cal P} (H_l)$.
(The critical points of
${\bar H}_l = - H_l$
generating 
$C_\ast ({\bar H}_l )$
are viewed here as critical points of 
$H_l$).

Observe that such a representation 
(\ref{eqn-Theta-as-a-chain-in-Morse-complex})
of
$\Theta_{H,J}$
is not unique since
different 
$\Lambda_\omega{\hbox{\rm -linear}}$ 
combinations of the critical points 
$z_j^l$
may represent the same homology class in
$H_\ast (M, \Lambda_\omega)$.

\smallskip
{\sl
Now we claim that for any
$\tau$
as above the number
$n (\hat{\gamma}, H, \tilde{J}_\tau )$
has to be non-zero and therefore
$\hat{\gamma}$
is durable.
}
\smallskip

Indeed, since the critical point
$z^l$
of  
${\bar H}_l$
is homologically essential for the rational singular homology class
Poincar{\'e}-dual to
$c_A$
one gets that the decomposition
(\ref{eqn-Theta-as-a-chain-in-Morse-complex})
must include a term
\[
\alpha_{j_1}^1\otimes\ldots\otimes \alpha_{j_{l-1}}^{l-1}\otimes 
{\hat{\gamma}}_l,
\]
such that
$c_i (\alpha_{j_i}^i)\neq 0$
for any
$i=1,\ldots ,l-1$.
Moreover, it follows from the definition of an 
$i{\hbox{\rm -friendly}}$
basis
$\{ \alpha_j^i\}$
(see
Definition~\ref{def-involved})
that such a term is unique. 
Comparing
(\ref{eqn-theta-chain-[M]})
and
(\ref{eqn-Theta-as-a-chain-in-Morse-complex})
one sees that the existence and the uniqueness of such a term in
(\ref{eqn-Theta-as-a-chain-in-Morse-complex})
imply that
$n (\hat{\gamma}, H, \tilde{J}_\tau )\neq 0$.
The claim is proven.

This finishes the proof of 
Proposition~\ref{prop-quantum-multipl-implies-nontriv-gw-number}
and 
Theorem~\ref{thm-main-time-indep-hams}.
\b
\smallskip

\vfil
\eject

\bigskip
\bigskip
\section{Proofs of 
Theorems~\ref{thm-two-small-hams}, 
\ref{thm-many-small-hams}, 
\ref{thm-3-slow-hams-two-dim-torus}} 
\label{pfs-of-corollaries-about-slow-hams}

In 
Theorems~\ref{thm-two-small-hams}
and
\ref{thm-many-small-hams}, 
we can, without loss of generality, consider only the special case
when all the slow 
Hamiltonians are Morse functions that have unique points of 
global maximum and minimum. Indeed,
this condition can be always achieved by
a sufficiently
$C^\infty{\hbox{\rm -small}}$
perturbation of the Hamiltonians as functions on
$M$.
Thus, since a sufficiently
$C^\infty{\hbox{\rm -small}}$
perturbation of a slow Hamiltonian is again slow,
the general case follows from the special case by continuity.

Now 
Theorem~\ref{thm-many-small-hams}
follows from 
Theorem~\ref{thm-main-time-indep-hams}
and
Example~\ref{exam-tuples-of-maxs-min-involved-in-multipl-ident}.

Theorem~\ref{thm-two-small-hams}
follows from
Theorem~\ref{thm-many-small-hams}
applied in the case
$l=2$
to the pairs
$H_1, -H_2$
and
$-H_1, H_2$.

To prove
Theorem~\ref{thm-3-slow-hams-two-dim-torus}
check, as in 
Example~\ref{exam-tuples-of-crit-pts-involved-in-multipl-ident},
that since the Hamiltonians are perfect Morse functions
$\hat{\gamma}\in {\cal P} (H)$
associated with
$z_1, z_2, z_3, A=0$
is involved in the identity
\[
c_1\cup c_2 = c_3,
\]
where
$c_1, c_2$
form the basis in
$H^1 ({\bf T}^2)$
dual to
$\alpha_1, \alpha_2$
(i.e. 
$c_i (\alpha_j) = \delta_{ij}$),
and
$c_3$
is the generator of
$H^2 ({\bf T}^2)$
which evaluates as 1 on the fundamental class.
Then
Theorem~\ref{thm-3-slow-hams-two-dim-torus}
follows 
Theorem~\ref{thm-main-time-indep-hams}.
\b
\smallskip

\vfil
\eject
\bigskip
\bigskip
\section{Proof of Theorem~\ref{thm-abw-ineqs} }
\label{sect-pf-thm-abw-ineqs}

Recall that the function 
$\Upsilon_l$
is on conjugacy classes in
$SU (n)$
was defined 
with respect to the Finsler metric defined by the Finsler norm
induced from
$T_\ast {\hbox{\it Ham}}\,  (Gr\, (r,n),\omega)$.
Then the first inequality follows by functoriality
directly from the definitions. 

Now we will prove the second inequality. Let us start with some 
technical observations.
As it was already mentioned, the Grassmannian
$Gr\, (r,n)$
can be viewed as the result of symplectic reduction
for the Hamiltonian action of
$U(r)$
(by multiplication from the right) on the space 
${\bf C}^{rn}$
of complex
$r\times n$
matrices.
Namely the value of the moment map 
$F: {\bf C}^{rn}\to {\bf R}$
of the action on an
$r\times n$
matrix
$A$
can be written as
\[
F (A) = B^\ast B/2 i
\]
and the Grassmannian 
$Gr\, (r,n)$
can be identified with the quotient
$S/ U(r)$,
where
$S = F^{-1} (Id /2i)\subset {\bf C}^{rn}$.

Given an
$\alpha = (\alpha_1, \ldots ,\alpha_n)\in {\mathfrak U}$
consider the Hamiltonian action 
$\{ g^\prime_{e^{2\pi i t}} \}$,
$e^{2\pi i t}\in S^1$,
of
$S^1$
on
${\bf C}^n$,
where an element
$e^{2\pi i t}\in S^1$
acts on a vector
$w= (w_1,\ldots ,w_n)\in {\bf C}^n$
by the formula
\[
e^{2\pi i t}: 
(w_1,\ldots, w_n)\to 
(e^{2\pi i t \alpha_1} w_1,\ldots ,  e^{2\pi i t \alpha_n} w_n).
\]
The action 
$\{ g^\prime_{e^{2\pi i t}} \}$
can be viewed as the Hamiltonian flow generated by the autonomous
Hamiltonian
$H^\prime_{\alpha} : {\bf C}^n\to {\bf R}$
defined by the formula
\[
H^\prime_{\alpha} (w) = -\alpha_1 |w_1 |^2 -\ldots -\alpha_n |w_n |^2. 
\]

The action
$\{ g^\prime_{e^{2\pi i t}} \}$
induces a Hamiltonian action
$\{ g_{e^{2\pi i t}} \}$,
$e^{2\pi i t}\in S^1$,
of
$S^1$
on
$Gr\, (r,n)$
which can be viewed as the Hamiltonian flow generated by some
{\it
normalized
}
Hamiltonian
$H_\alpha : Gr\, (r,n)\to {\bf R}$.

To study the properties of
$H_\alpha$
consider first the Hamiltonian action 
$\{ g^{\prime\prime}_{e^{2\pi i t}} \}$,
$e^{2\pi i t}\in S^1$,
of
$S^1$
on
${\bf C}^{rn}= ({\bf C}^n )^r$
obtained as the direct product of the actions 
$\{ g^\prime_{e^{2\pi i t}} \}$
on the factors
${\bf C}^n$
of
$({\bf C}^n )^r$.
The action 
$\{ g^{\prime\prime}_{e^{2\pi i t}} \}$,
can be viewed as the Hamiltonian flow generated by the Hamiltonian
$H^{\prime\prime}_\alpha : ({\bf C}^n )^r\to {\bf R}$
which is the direct sum of the Hamiltonians
$H^\prime_\alpha$
on the 
$r$
factors of
$({\bf C}^n )^r$.

The Hamiltonian
$H^{\prime\prime}_\alpha$
is constant along the fibers of the action of
$U (r)$
on
$S$
and thus descends to a function on
$Gr\, (r,n)$.
Since
$\alpha\in {\mathfrak U}$,
one has that
$\alpha_1 +\ldots +\alpha_n = 0$
and from here one easily deduces that the integral of
$H^{\prime\prime}_\alpha$
over
$S$,
and hence the integral of
$H_\alpha$
over
$Gr\, (r,n)$
are zero.
Thus the Hamiltonian
$H^{\prime\prime}_\alpha$
restricted on
$S$
descends exactly to the 
{\it
normalized
}
Hamiltonian 
$H_\alpha$
on
$Gr\, (r,n)$.

Assume now that
$\alpha_{j_1}\neq \alpha_{j_2}$
for
$j_1\neq j_2$.
This implies that
$H_\alpha$
is a Morse function on
$Gr\, (r,n)$.
From the definition of 
${\mathfrak U}$
one sees that
$H_\alpha$
also satisfies the conditions 
$(B)$
and
$(C)$
of 
Definition~\ref{def-slow-hams}.
Thus
$H_\alpha$
is a slow Morse Hamiltonian.

The considerations above allow us to compute easily the critical values of
$H_\alpha$.
At a critical point
$z_I$
corresponding to a subset
$I= \{ i_1,\ldots ,i_r\} \subset \{ 1,\ldots ,n\}$
the critical value of
$H_\alpha$
is equal to
\[
H_\alpha (z_I) = - \sum_{j\in I} \alpha_j.
\]
Hence
\begin{equation}
\label{eqn-action-Grassm-case}
{\cal A}_\alpha  ([z_I, f_d]) = \sum_{j\in I} \alpha_j - d,
\end{equation}
where
$[z_I, f_d]\in {\cal P} (H_\alpha )$
is formed by a constant path 
$z_I$
and a two-dimen\-sion\-al sphere
$f_d$
attached to
$z_I$
and representing the class
$d\in {\bf Z} \cong \Pi$.

The critical points of
$H_\alpha$
are the fixed points of the action
$\{ g_{e^{2\pi i t}} \}$
and therefore they are in one-to-one correspondence with
the invariant complex
$r{\hbox{\rm -dimensional}}$
subspaces of the matrix
$diag (e^{2\pi i\alpha_1 t},\ldots ,e^{2\pi i\alpha_n t})$
or, equivalently, with the subsets
$I= \{ i_1,\ldots ,i_r\}\subset \{ 1,\ldots , n \}$.
Thus the total number of the critical points of
$H_\alpha$
is equal to the total sum of the Betti numbers of
$Gr\, (r,n)$
and therefore
$H_\alpha$
is a 
{\it
perfect
}
Morse function on
$Gr\, (r,n)$.

Let us now consider the isomorphism between the Morse homology of
$H_\alpha$
and the singular homology of
$Gr\, (r,n)$
and let us show that a critical point
$z_I$
corresponds to the homology class
$\sigma_I$.
Assume without loss of generality that the flag 
$\{0\} = F_0\subset\ldots \ldots\subset F_n = {\bf C}^n$
used to define
the homology classes
$\sigma_I$
is the standard flag spanned by the basic vectors 
${\bf e}_1,\ldots , {\bf e}_r$
associated with the
coordinates
$w_1,\ldots, w_l$
on
${\bf C}^n$:
\[
F_k = {\rm span}\, ({\bf e}_1,\ldots , {\bf e}_k),
\]
$k=1,\ldots, n$.
Given a subset
$I= \{ i_1,\ldots ,i_r\}\subset \{ 1,\ldots , n \}$
consider the 
$r{\hbox{\rm -dimensional}}$
subspace
$V_I\subset {\bf C}^n$
generated by the 
$r$
basic vectors 
${\bf e}_{i_1},\ldots , {\bf e}_{i_r}$
and the
$(n-r){\hbox{\rm -dimensional}}$
subspace 
$V_I^{\perp}$
spanned by the rest of the basic vectors
${\bf e}_j$.
A neighborhood of
$V_I$
in
$Gr\, (r,n)$
can be identified with the space
${\bf C}^{r(n-r)}$
of 
$r{\hbox{\rm -tuples}}$
$({\bf v}_1,\ldots ,{\bf v}_r)$,
${\bf v}_j\in V_I^\perp$,
$j=1,\ldots ,r$:
namely, to a tuple
${\bf v} = ({\bf v}_1,\ldots ,{\bf v}_r)$
one associates the point
$x_{\bf v} \in Gr\, (r,n)$
corresponding to the 
$r{\hbox{\rm -dimensional}}$
subspace of
${\bf C}^n$
spanned by
$({\bf e}_1 + {\bf v}_1,\ldots , {\bf e}_r + {\bf v}_r )$.
The function
$f$
on a neighborhood of
$x_0 = V_I$
in
$Gr\, (r,n)$
defined by the formula
\[
f (x_{\bf v}) = H^{\prime\prime}_\alpha 
( {\bf e}_1 + {\bf v}_1,\ldots , {\bf e}_r + {\bf v}_r )
\]
has a Morse critical point at
$x_0$
and, under our choices of coordinates, the quadratic parts of
$f$
and
$H_\alpha$
at
$x_0 = V_I$
coincide.
This can be seen if one applies the Gram-Schmidt procedure to turn the
basis 
$( {\bf e}_1 + {\bf v}_1,\ldots , {\bf e}_r + {\bf v}_r )$
of the subspace
$x_{\bf v}$
into a unitary basis which can be viewed as an element of
$S$.
It allows us to check that the unstable manifold of
$H_\alpha$
at
$x_0 = V_I$,
taken with respect to the natural Riemannian metric defined by the complex 
and the symplectic structures on the Grassmannian,
coincides with
$W_I$.
Using this fact one gets
that under the identification 
of the (integral) Morse homology of
$-H_{\zeta^j}$
with the singular homology of
$Gr\, (r,n)$ 
the homology class
$z_I\in H_\ast (C_\ast (-H_{\zeta^j} ))$
represented by
$z_I$
viewed as a critical point of
$-H_{\zeta^j}$
gets identified with the class
$\sigma_{\ast I}$,
where
$\ast I = \{ n+1- i_1,\ldots , n+1-i_r\} \subset \{ 1,\ldots ,n\}$
and 
$\sigma_{\ast I}$
is the homology class defined by
$\ast I$
which is in fact Poincar{\'e}-dual to
$\sigma_I$.

Now let us finish the proof of the second inequality of the theorem.
We prove it first for a
{\it
good
}
$\zeta$,
i.e. for a
$\zeta\in {\mathfrak U}^l$
such that for each
$\zeta^i = (\zeta_1^i,\ldots ,\zeta_n^i)\in {\mathfrak U}$,
$1\leq i\leq l$,
one has
$\zeta_{j_1}^i \neq \zeta_{j_2}^i$
for any
$j_1\neq j_2$.
Such a 
$\zeta$
gives rise to slow Morse Hamiltonians
$H_\zeta = (H_{\zeta^1},\ldots , H_{\zeta^l} )$.
Using our computations of the critical values of a slow Morse Hamiltonian
$H_{\zeta^i}$,
$i=1,\ldots ,l$,
and the identification of its critical points with the homology
classes of the Grassmannian as above one applies
Corollary~\ref{cor-thm-main-time-indep-hams-perfect-hams}
to obtain the second inequality of the theorem in the case when
$\zeta$
is good.
Finally observe that the set of good
$\zeta$
is open and dense in
${\mathfrak U}^l$.
Therefore by continuity one gets the second inequality in the general case.
\b
\smallskip

\vfil
\eject

\bigskip
\bigskip
\section{Proofs of the results concerning pseudo-holo\-mor\-phic curves}
\label{pfs-results-concern-psh-curves}

We will briefly outline the proofs of the results concerning the 
moduli spaces of pseudo-holomorphic curves.

\bigskip
\subsection{The proof of 
Proposition~\ref{prop-moduli-space-generically-smooth}
}
\label{pf-prop-moduli-space-generically-smooth}

The basic scheme of the proof is fairly standard for the Floer theory
(see
\cite{Flo},
\cite{Ho-Sa}
or
\cite{Sal}).
We outline it for our particular case following M.Schwarz 
-- the details can be found in
\cite{Sch-PhD},
\cite{Ho-Sa}.

\bigskip
\noindent
{\bf Preliminaries.}

\smallskip
One defines 
$C^\infty_{\gamma} (\Sigma, M)$
as the space of smooth maps
$F:\Sigma\to M$
such that over the end
$\Sigma_i$,
$i=1,\ldots ,l$,
one has
\[
F\circ\Phi_i (\frac{s}{\sqrt{1-s^2}}, t) = \phi_i (s,t)
\]
for some 
$C^\infty{\hbox{\rm -function}}$
$\phi_i$
on the 
{\it closed}
cylinder
$[0,1]\times S^1$,
where
$\phi_i (1,\cdot) = \gamma_i$.
The space
$C^\infty_{\gamma} (\Sigma, M)$
can be completed to a separable infinite-dimensional Banach manifold 
$H^{1,p}_{\gamma} (\Sigma, M)$
modeled on the Sobolev space
$H^{1,p} (\Sigma, {\bf R}^{2n})$
with
$p>2$.

Consider now the special case when
$\Sigma = {\bf R}\times S^1$
is the standard cylinder.
Let
$h: S^1\times M\to {\bf R}$
be a Hamiltonian function.
We will say that
$h$
is
{\it
non-degenerate,
}
if
all the contractible 1-periodic orbits of the 
Hamiltonian flow of
$h$
are non-degenerate
(i.e. 1 is not a Floquet multiplier for any of the orbits).
A generic Hamiltonian is non-degenerate.
Now let
$J$
be an
$\omega{\hbox{\rm -compatible}}$
almost complex structure on 
$M$.
Let
$\gamma^\prime, \gamma^{\prime\prime}$
be contractible 1-periodic orbits of the Hamiltonian flow of
$h$.
Define the operator
\[
\partial_{J,h} = \partial_s + J\partial_t -\nabla h
\]
on  
$H^{1,p}_{\gamma^\prime, \gamma^{\prime\prime}} 
({\bf R}\times S^1, M)$,
where gradient is taken with respect to the Riemannian metric 
$g (\cdot ,\cdot ) =\omega (\cdot , J\cdot )$
on
$M$.
Consider the space
${\cal M} (\gamma^\prime, \gamma^{\prime\prime})\subset 
H^{1,p}_{\gamma^\prime, \gamma^{\prime\prime}} ({\bf R}\times S^1, M)$
of solutions 
$u$
of the equation
\[
\partial_{J,h} u = 0.
\]
We say that the pair
$(h, J)$
is 
{\it
weakly regular,
}
if
$h$
is non-degenerate and the linearization of
$\partial_{J,h}$
at any point
$u\in {\cal M} (\gamma^\prime, \gamma^{\prime\prime})$
is onto.
A generic (with respect to the
$C^\infty{\hbox{\rm -topology}}$)
pair 
$(h, J)$
is weakly regular
\cite{Sal-Ze}.
Moreover, if
$h$
is a slow Morse Hamiltonian which is a Morse-Smale function with respect to 
the metric 
$g (\cdot ,\cdot ) =\omega (\cdot , J\cdot )$
then the pair
$(h, J)$
is weakly regular
\cite{Sal-Ze}.
Given the Hamiltonians
$H = (H_1,\ldots ,H_l)$
and the almost complex structure
$\tilde{J}$
as above we say that the pair
$(H, \tilde{J})$
is
{\it
weakly regular
}
if for each
$i=1,\ldots ,l$
the pair
$(H_i, J (\tilde{J}) )$
is weakly regular.

\bigskip
\noindent
{\bf Smoothness of
${\cal M} (\hat{\gamma}, H, \tilde{J})$.
}

The operator
${\bar{\partial}}_{\tilde{J}}$
defined on
$C^\infty_{\gamma} (\Sigma, M)$
can be extended to
$H^{1,p}_{\hat{\gamma}} (\Sigma, M)$.
Consider the solutions
$u\in H^{1,p}_{\hat{\gamma}} (\Sigma, M)$
of the equation
${\bar{\partial}}_{\tilde{J}} u = 0$.
If
$(H, \tilde{J})$
is weakly regular then using local 
elliptic regularity and Sobolev embedding theorems
one shows that any such solution has to lie inside
$C^\infty_{\gamma} (\Sigma, M)$.

Let us view elements of
${\cal M} (\hat{\gamma}, H, \tilde{J})$
as maps 
$\Sigma\to M$.
Then
${\cal M} (\hat{\gamma}, H, \tilde{J})$
can be viewed as the zero level set of the operator
${\bar{\partial}}_{\tilde{J}}$
defined on 
$H^{1,p}_{\hat{\gamma}} (\Sigma, M)$.

Using the Fredholm theory one shows that 
if
$(H,\tilde{J})$
is weakly regular then the linearization of the operator
${\bar{\partial}}_{\tilde{J}}$
on
$H^{1,p}_{\hat{\gamma}} (\Sigma, M)$
at any 
$u\in {\cal M} (\hat{\gamma}, H, \tilde{J})$
is Fredholm for any
$p\geq 2$.
For 
$p=2$
the proof of this fact is fairly easy
(see e.g.
\cite{Flo},
\cite{Sch-PhD}). 
For
$p>2$
the result is deduced in
\cite{Sch-PhD}
from the case
$p=2$
by means of a local 
$L^p{\hbox{\rm -estimate}}$ 
for the Cauchy-Riemann operator
(see e.g.
\cite{Sch-PhD},
\cite{Sal}).

To compute the index of the elliptic operator
${\bar{\partial}}_{\tilde{J}}$
one uses the index additivity with respect to gluing of 
trivial bundles over different surfaces along with the operators
${\bar{\partial}}_{\tilde{J}}$
on those bundles.
Namely, one computes directly the index for the operator 
${\bar{\partial}}_{\tilde{J}}$
on the bundle over a plane (an open disc with the cylindrical end).
Then one caps off
$\Sigma$
with the discs (gluing together the bundles and the operators 
over the cylindrical ends) and reduces the problem to the case of a
${\bar{\partial}}{\hbox{\rm -operator}}$
over a closed Riemann surface, where the answer is given by the
the Riemann-Roch formula.
Finally one gets that the index of the operator
${\bar{\partial}}_{\tilde{J}}$
on the bundle over
$\Sigma$
is equal to
$2n -\mu (\hat{\gamma} )$.

Thus the implicit function theorem for Banach manifolds implies that if
the linearization of
${\bar{\partial}}_{\tilde{J}}$
is onto at any
$u\in {\cal M} (\hat{\gamma}, H, \tilde{J})$
then
${\cal M} (\hat{\gamma}, H, \tilde{J})$
is a smooth manifold of dimension
$2n -\mu (\hat{\gamma})$.
In such a case we say that a (weakly regular) pair
$(H, \tilde{J})$
is
{\it
regular.
}

\bigskip
\noindent
{\bf Transversality.}

A transversality result based on
an infinite-dimensional version of the Sard 
theorem shows that
regular pairs form a residual set in the set of all weakly regular 
pairs. 
Here is the first place where we have to do something to put the discussion
in the framework of
\cite{Sch-PhD} 
(recall that our class of almost complex structures is larger than in
\cite{Sch-PhD}) -- 
all the previous statements were literally covered by the results in
\cite{Sch-PhD}. 
We can assume without loss of generality 
that the symplectic form
$\Omega_{\tilde{\omega}, \tau }$
on
$\Sigma\times M$,
with which our almost complex structures 
$\tilde{J}$
are compatible, coincides with the split form
$\pi^\ast\Omega\oplus \omega$ 
outside of a compact set contained in the complement
$\Sigma^0 = \Sigma\setminus\bigcup_{i=1}^l \Sigma_i$.
In 
\cite{Sch-PhD}
the strategy was to perturb the extension of the Hamiltonians
$H_1,\ldots , H_l$
over 
$\pi^{-1} (\Sigma^0)$
that had been defined by means of the cut-off function
$\beta$
in order to perturb the term corresponding to the gradients of the
Hamiltonians in the
$\bar{\partial}{\hbox{\rm -operator}}$.
Then one considered in
\cite{Sch-PhD}
variations of the almost complex structures 
$J_x$
on the fibers
$\pi^{-1} (x)$,
$x\in\Sigma^0$,
to get the necessary result.
In our case, due to our definition of 
${\cal T} (H)$,
the perturbation of the Hamiltonian term in the
$\bar{\partial}{\hbox{\rm -operator}}$
can be viewed as the result of a perturbation of the
almost complex structure
$\tilde{J}\in {\cal T} (H)$.
Then one follows the proof of
Theorem 4.2.20 in
\cite{Sch-PhD} 
and obtains the following result: 
given a weakly regular pair 
$(H, \tilde{J})$
one can always find an open set  
$K\subset \Sigma$
with compact closure so that for a generic and arbitrarily
$C^\infty{\hbox{\rm -small}}$
perturbation 
$\tilde{J}_1\in {\cal T} (H)$
of
$\tilde{J}$
over
$K$
the pair
$(H, \tilde{J}_1 )$
is regular.
Observe also that a sufficiently small perturbation 
$\tilde{J}_1\in {\cal T} (H)$
of
$\tilde{J}$
remains compatible with the symplectic form
$\Omega_{\tilde{\omega}, \tau }$.

\bigskip
\noindent
{\bf Compactness.}

As before let
$\pi: \Sigma\times M\to\Sigma$,
$pr_M :T_\ast (\Sigma\times M)\to T_\ast M$
be the natural projections and
let 
$J_x$
be the restriction of
$\tilde{J}$
on
$\pi^{-1} (x)$.
Let the norm
$\|\cdot\|_x$
on the tangent bundle of the fiber
$\pi^{-1} (x)$
be defined by the metric
$\omega (\cdot, J_x \cdot)$.

Given a map
$u:\Sigma\to \Sigma\times M$
from
${\cal M} (\hat{\gamma}, H, \tilde{J})$
and a conformal chart
$f: U\to \Sigma$,
$U\in {\bf C}$,
with coordinates 
$s, t$
set
\[
u_s = pr_M\circ \displaystyle{\frac{\partial (u\circ f)}{\partial s}},
u_t = pr_M\circ \displaystyle{\frac{\partial (u\circ f)}{\partial t}}
\]
and define the 2-form 
$\| du\|\Omega$
on
$f (U)$,
where
$\Omega$
is the volume form on
$\Sigma$,
as
\[
\| du\|\Omega (x) = f^\ast \{ ( {\| u_s (x) \|}_x + {\| u_t (x) \|}_x ) 
ds\wedge dt \}.
\] 
One can check that the form
$\| du\|\Omega$
defined in such a way is a correctly defined 2-form on the whole
$\Sigma$.
Define the 
{\it
energy
}
of
$u$
as
\[
E (u) = \int_\Sigma \| du\| \Omega.
\]

Compactness of the space
${\cal M} (\hat{\gamma}, H, \tilde{J})$
can be shown by the following standard argument
(see
\cite{Ho-Sa},
\cite{HZ},
\cite{Sch-PhD}).
First of all one checks that there is an
{\it 
a priori 
}
uniform bound on the energies of all 
$u\in {\cal M} (\hat{\gamma}, H, \tilde{J})$.
(see
\cite{Sch-PhD}, 
also see
\cite{Akv-Sal}, Lemma 5.2).
Now there are only two phenomena that might obstruct the compactness of
${\cal M} (\hat{\gamma}, H, \tilde{J})$:
bubbling-off of 
$\tilde{J}{\hbox{\rm -holomorphic}}$
spheres or convergence to a ``broken trajectory''
(see e.g.
\cite{Flo},
\cite{Ho-Sa},
\cite{HZ},
\cite{Sal}
or
\cite{Sch-PhD}).
The latter obstacle does not occur in our case since we consider only
zero-dimensional moduli spaces
${\cal M} (\hat{\gamma}, H, \tilde{J})$
(otherwise one has to use Floer's gluing techniques -- see
\cite{Flo}, 
also see
\cite{Sch-book}).
As far as the bubbling-off is concerned one quickly observes that 
because of our definition of
${\cal T} (H)\supset \tilde{J}$
the projection of such a bubble on
$\Sigma$
has to be holomorphic and thus the maximum principle dictates that
the bubble has to lie inside a fiber of
$\Sigma\times M\to\Sigma$.
Then one uses the fact that
$(M,\omega)$
is strongly semi-positive and therefore
(see
\cite{Ho-Sa},
\cite{Se})
a generic 2-parametric family of 
$\omega{\hbox{\rm -compatible}}$
almost complex structures on
$(M,\omega)$
does not contain an almost complex structure 
$J$
which admits a
$J{\hbox{\rm -holomorphic}}$
sphere of negative Chern number.
Such a 2-parametric family arises from
$\tilde{J}$
since the restriction of
$\tilde{J}$
may vary with the base point in the two-dimensional surface
$\Sigma$.
Also because of the transversality reasons one can always assume without loss
of generality that the 1-periodic trajectories of the Hamiltonian flows of
$H_1, \ldots, H_l$
do not intersect the
$J(\tilde{J}){\hbox{\rm -holomorphic}}$
spheres with Chern number 1 in
$M$.
Then one proceeds as in
\cite{Ho-Sa}
(cf.
\cite{Sch-cup})
and shows that for such a 
$\tilde{J}$
in the absence of pseudo-holomorphic curves with negative Chern numbers 
the bubbling-off in the fibers cannot occur, which shows the compactness of
${\cal M} (\hat{\gamma}, H, \tilde{J})$.

\bigskip
\noindent
{\bf Orientation.}

The orientation can be obtained by the methods from
\cite{Flo-Ho}.

\bigskip
\noindent
This finishes the proof of 
Proposition~\ref{prop-moduli-space-generically-smooth}.
\b
\bigskip

\bigskip
\subsection{The proof of 
Theorem~\ref{thm-main-pss}}
\label{pf-thm-main-pss}

Pick
a connection
$\nabla\in {\cal L}_{\{ c\} }^{fit} ([\varphi_H])$,
the corresponding 2-form
$\tilde{\omega}_\nabla = \tilde{\omega}$,
such that
$\tau_0 < {\hbox{\it size}}\, (\tilde{\omega} )$,
and a weak coupling deformation
$\{\Omega_{\tilde{\omega}, \tau^\prime }\}$,
$0< \tau^\prime < \tau$.
Observe that for 
$\tau^\prime$
close to zero almost complex structures from
${\cal T}_0 (H)$
are compatible with
$\Omega_{\tilde{\omega}, \tau^\prime }$.

Now the transversality and compactness
results from the proof of
Proposition~\ref{prop-moduli-space-generically-smooth}
in
Section~\ref{pf-prop-moduli-space-generically-smooth}
also hold for 1-parametric families of almost complex structures from
${\cal T} (H)$.
(The compactness result relies on the fact that
a generic 2-parametric family of 
$\omega{\hbox{\rm -compatible}}$
almost complex structures on
$(M,\omega)$
does not contain an almost complex structure 
$J$
which admits a
$J{\hbox{\rm -holomorphic}}$
sphere of negative Chern number).
This allows us to choose some small
$\epsilon > 0$
and a family 
$\{ \tilde{J}_{\tau^\prime}\}$,
$\epsilon\leq\tau^\prime \leq \tau$,
so that: 

\medskip
\noindent
$\bullet$
Each
$\tilde{J}_{\tau^\prime}$ 
belongs to 
${\cal T}_J (H)$
and is compatible with
$\Omega_{\tilde{\omega}, \tau^\prime }$.

\smallskip
\noindent
$\bullet$
For all 
$\tau^\prime$
sufficiently close to
$\epsilon$
one has
$\tilde{J}_{\tau^\prime}\in {\cal T}_0 (H)$.

\smallskip
\noindent
$\bullet$
For all 
$\tau^\prime$
sufficiently close to
$\epsilon$
and to
$\tau^\prime$
the pair
$(H,\tilde{J}_{\tau^\prime})$ 
is regular.

\smallskip
\noindent
$\bullet$
For any 
$\hat{\gamma}$
with 
$\mu (\hat{\gamma}) = 2n$
the union 
\[
\tilde{\cal M} (\hat{\gamma}, H, \tilde{J}_{\tau^\prime}) 
=\bigcup_{\varepsilon\leq \tau^\prime \leq \tau}
{\cal M} (\hat{\gamma}, H, \tilde{J}_{\tau^\prime})
\]
forms a compact (oriented) 1-dimensional manifold with boundary 
(where
${\cal M} (\hat{\gamma}, H, \tilde{J}_\varepsilon )$
and
${\cal M} (\hat{\gamma}, H, \tilde{J}_\tau)$
are among the connected components of the boundary).

\smallskip
\noindent
$\bullet$
For any 
$\hat{\gamma}$
with 
$\mu (\hat{\gamma}) = 2n+1$
the union 
\[
\tilde{\cal M} (\hat{\gamma}, H, \tilde{J}_{\tau^\prime}) 
=\bigcup_{\varepsilon\leq \tau^\prime \leq \tau}
{\cal M} (\hat{\gamma}, H, \tilde{J}_{\tau^\prime})
\]
is a compact (oriented) 0-dimensional manifold.
\medskip

Now for
$\tau^\prime$
close to
$\epsilon$,
when
$\tilde{J}_\tau\in {\cal T}_0 (H)$,
the theorem is literally the result from
\cite{PSS}.
To extend it to other
$\tau^\prime$
one uses the standard homotopy argument from the Floer theory as in
\cite{Sch-PhD} (Ch. 5.2). 
Namely, one 
constructs a chain endomorphism of
$CF_\ast (H_1, J)\otimes\ldots\otimes  CF_\ast (H_l, J)$
by counting points in zero-dimensional spaces 
$\tilde{\cal M} (\hat{\gamma}, H, \tilde{J}_{\tau^\prime})$.
This endomorphism maps the chain
$\theta_{\Sigma, H,\tilde{J}_\epsilon }$
into the chain
$\theta_{\Sigma, H,\tilde{J}_\tau }$.
Then one constructs a chain homotopy of the endomorphism to the identity
by means of the 1-dimensional spaces 
$\tilde{\cal M} (\hat{\gamma}, H, \tilde{J}_{\tau^\prime})$
thus showing that the homology classes of the cycles realized by
$\theta_{\Sigma, H,\tilde{J}_\epsilon }$
and
$\theta_{\Sigma, H,\tilde{J}_\tau }$
coincide.
\b
\bigskip

\bigskip
\subsection{The proof of 
Proposition~\ref{prop-isom-betw-Morse-Floer-sing-hom-for-slow-Hams}}
\label{pf-prop-isoms-betw-Morse-Floer-complex-commute}

We use the same kind argument that was used in 
\cite{PSS}
in the proof of the statement that
${\cal QF}_h$
is an isomorphism between quantum
and Floer (co)homology. 
For a slow Morse Hamiltonian
$h$
one can choose a
{\it
time-independent
}
almost complex structure
$J$
on
$M$
so that the pair
$(h, J)$
is regular.

To define the map
${\cal QF}_h$
(see
\cite{PSS},
the proof of Theorem 4.1)
one first identifies the quantum homology with the homology of
the Morse complex of
$h$
with coefficients in
$\Lambda_\omega$.
Since the Morse function
$h$ 
is slow its Morse complex can be canonically identified with its Floer 
complex
$CF_\ast (h,J)$.
Then one constructs
${\cal QF}_h$
as the map in homology induced by a chain endomorphism of that complex.
(For simplicity we denote by
${\cal QF}_h$
the induced maps both in homology and cohomology).
Thus we need to prove that this chain endomorphism induces the identity
map in the homology of the Morse complex. We will show that it is actually
identity already on the level of chains.

Indeed, recall how  
${\cal QF}_h$
is defined as a chain endomorphism. One considers the moduli space 
${\cal M}_{\bf C} (y, h, J)$
of ``spiked disks'', i.e. the space defined as in
Definition~\ref{def-modul-space}
with
$l=1$,
$\Sigma = {\bf C}$,
where
$y$
is a critical point of
$h$.
Then one picks critical points
$x,y$
of
$h$
in such a way that the space of curves from
${\cal M}_{\bf C} (y, h, J)$
with
$u (0)$
lying on the unstable manifold of a critical point 
$x$
is zero-dimensional. Counting such curves in various spherical non-torsion
homology classes
$A$,
such that
${\rm ind}\, z -  {\rm ind}\, x = 2 c_1 (A)$,
one defines the coefficients 
$\langle x, y, A\rangle$
that determine the map
${\cal QF}_h$.

Consider the deformation
$\tau h$, 
$0\leq\tau\leq 1$.
The number 
$\langle x, y, A\rangle$,
defined by means of the Hamiltonian 
$\tau h$,
does not change as
$\tau$
goes from one to zero.
But when
$\tau$
is zero the number  
$\langle x, y, A\rangle$
counts 
$J{\hbox{\rm -holomoprhic}}$ 
curves in the homology class 
$A$.
Therefore, since
$(M,\omega )$
is strongly semi-positive,
$c_1 (A)$
has to be non-negative for any 
$A$
such that
$\langle x, y, A\rangle \neq 0$.
 
Now, if one uses a time-independent
$J$
then, since
$h$
is also time-in\-depen\-dent, any curve from 
${\cal M}_{\bf C} (y, h, J)$
with
$u (0)$
lying on the unstable manifold of a critical point 
$x$
can be ``rotated'' by changing the time parameter
and in this way one gets a one-dimensional family of solutions. This leads to
a contradiction with the zero dimension of
${\cal M}_{\bf C} (y, h, J)$, 
unless the original curve which we rotated was time-independent and hence
was a part of the anti-gradient trajectory of
$h$
going from
$x$
to
$y$.
But since
$c_1 (A)$
is always non-negative,
${\rm ind}\, z -  {\rm ind}\, x \geq 0$
and therefore we must have
$x = y$, 
$A =0$,
which shows that
${\cal QF}_h$
is identity already on the level of chains.
\b
\bigskip

\vfil
\eject

\bigskip
\bigskip
\section{Proof of 
Proposition~\ref{prop-upsilon-via-systems-of-paths}}
\label{sect-pf-prop-upsilon-via-systems-of-paths}

Given a tuple
${\cal C} = ({\cal C}_1 ,\ldots , {\cal C}_l )$
of conjugacy classes in
$G$
and a tuple 
$\phi = (\phi_1,\ldots ,\phi_l )$
of elements from
$G$
we will write
$\phi\in {\cal C}$
if
$\phi_i\in {\cal C}_i$,
$i=1,\ldots ,l$.
Denote by
$\Delta$
the set of all tuples
$f = (f_1,\ldots, f_l)$
of elements from
$G$
such that  
$f_1\cdot\ldots\cdot f_l = Id$.

Consider systems 
of paths 
$a = (a_1,\ldots , a_l)\in {\cal G} ({\cal C})$ 
containing
$l-1$
constant paths 
$a_i$,
$i=2,\ldots , l$,
identically equal to
$\phi_i\in {\cal C}_i$.
Set
$\phi_1 = a_1 (1)\in {\cal C}_1$.
Then
\[
{\hbox{\it length}}\, (a) = {\hbox{\it length}}\, (a_1)\geq
\rho\, (Id, a_1 (1)\cdot a_1^{-1} (0)) = \rho\, (Id, \prod_{i=1}^l \phi_i).
\]
Thus we get that
\[
\Upsilon_{l}\, ({\cal C}) = \inf_{\phi\in {\cal C}} 
\rho\, (Id, \prod_{i=1}^l \phi_i)\geq 
\inf_{a\in {\cal G} ({\cal C})} {\it length}\, (a).
\]

Let us prove the opposite inequality. Because of the bi-invariance of
the pseudo-metric and the elementary property
$\Upsilon_{l}\, ({\cal C}) = \Upsilon_{l}\, ({\cal C}^{-1})$ 
it suffices to show that
\begin{equation}
\label{eqn-dist-cal-N-upsilon}
\inf_{f\in\Delta,\psi\in {\cal C}} 
\sum_{i=1}^l \rho\, (Id, f_i\psi_i^{-1} )\geq
\inf_{\phi\in {\cal C}} 
\rho\, (Id, \prod_{i=1}^l \phi_i^{-1}).
\end{equation}
\noindent
Using the triangular inequality one gets:
\[
\inf_{f\in\Delta,\psi^{-1}\in {\cal C}} \sum_{i=1}^l 
\rho\, (Id, f_i\psi_i^{-1})\geq
\inf_{f\in\Delta,\psi\in {\cal C}}  
\rho\, (Id, \prod_{i=1}^l f_i\psi_i^{-1} )
\]
Now set
$F_i =\prod_{j=1}^i f_j$,
$i=1,\ldots, l$,
use the identity
$F_l = \prod_{i=1}^l f_i = Id$
and observe that
\[
\prod_{i=1}^l f_i\psi_i^{-1} = 
\bigg(
\prod_{i=1}^l F_i \psi_i^{-1} F_i^{-1}
\bigg)
\cdot F_l
= \phi_1^{-1}\cdot\ldots\cdot\phi_l^{-l}
\]
for some 
$\phi = (\phi_1,\ldots,\phi_l )\in {\cal C}$.
This implies
(\ref{eqn-dist-cal-N-upsilon})
and the proposition is proven.
\b
\smallskip

\vfil
\eject

\vfil
\eject

\bigskip
\bigskip
\section{Proof of 
Theorem~\ref{thm-k-area-distance}}
\label{sect-proof-thm-k-area-distance}

\bigskip
\subsection{The case of a cylinder: 
$l=2$}
\label{sect-proof-thm-k-area-distance-cyl}

In the case when
$\Sigma$
is a cylinder
the proof basically imitates the similar proofs from
\cite{Pol1},
\cite{Pol3}.
\smallskip

\bigskip
\subsubsection{Coarse length}
\label{subsubsect-def-coarse-length}

Let
$\gamma : [0,1]\to G$
be a smooth path in
$G$.
Let
$\|\cdot \|$
be the Finsler norm on the tangent bundle of
$G$
defining our bi-invariant Finsler pseudo-metric on the group. 
Besides the usual length which is defined as
\[
{\hbox{\it length}}\, (\gamma) = \int_0^1 \|\frac{d\gamma}{ds}\|ds  
\]
one can define a quantity
${\hbox{\it coarse-length}}\, (\gamma)$
as
\[
{\hbox{\it coarse-length}}\, (\gamma) 
= \max_{s\in [0,1]} \|{d\gamma}/{ds}\|.
\]
Obviously one always has
\[
{\hbox{\it length}}\, (\gamma)
\leq
{\hbox{\it coarse-length}}\, (\gamma).
\]
Observe that 
${\hbox{\it coarse-length}}\, (\gamma)$,
unlike
${\hbox{\it length}}\, (\gamma)$,
essentially depends on the parameterization of the path
$\gamma: [0,1]\to G$.
However, by reparameterizing 
$\gamma$
one can always make its coarse length equal to its length
(see e.g.
\cite{Pol4}).
\smallskip

\bigskip
\subsubsection{Preliminaries}
\label{subsubsect-cylinder-prelim}

Without loss of generality we can assume that
$\Sigma$
is a standard cylinder
$[0,1]\times S^1$
with the coordinates
$(s,t)$,
$0\leq s\leq 1$,
$0\leq t\leq 1\, ({\hbox{\rm mod}}\, 1)$,
and equipped with the area form 
$\Omega = ds\wedge dt$,
so that
$\int_\Sigma \Omega = 1$.

Set
${\cal C} = ({\cal C}_1, {\cal C}_2 )$.
Recall from
Definition~\ref{def-cal-L}
that
${\cal L} ({\cal C})$
is the set of connections on
$P=\Sigma\times F\to\Sigma$
which are flat near the boundary and whose holonomies over the closed paths
${\{ t\to 0\times t \}}_{0\leq t\leq 1}$
and
${\{ t\to 1\times t \}}_{0\leq t\leq 1}$
belong, respectively, to the conjugacy classes
${\cal C}_1^{-1}$
and
${\cal C}_2$.
We denote by
${\cal G}^\prime ({\cal C})$
the set of all smooth paths
$a: [0,1]\to G$
such that
$a (0) \in {\cal C}_1^{-1}$,
$a (1) \in {\cal C}_2$.
By
${\cal G}_{\{  a\}  }^\prime ({\cal C})$
we denote the connected component of
${\cal G}^\prime ({\cal C})$
corresponding to the homotopy class
$\{  a\} $
of paths connecting
${\cal C}_1^{-1}$
and
${\cal C}_2$.

In our case, when
$l=2$,
the sets 
${\cal G} ({\cal C})$
and
${\cal G}^\prime ({\cal C})$
are closely related.
Indeed, an element of
${\cal G} ({\cal C}_1 , {\cal C}_2 )$
is a pair of paths
$a_1 , a_2: [0,1 ]\to G$
such that
$a_1^{-1} (0) = a_2 (0)$,
$a_1 (1) = \phi_1 \in {\cal C}_1 $,
$a_2 (1) = \phi_2 \in {\cal C}_2 $.
Joining the curves 
$a_1^{-1}$
and
$a_2$
in
$G$
at their common point
$a_1^{-1} (0) = a_2 (0)$
one obtains a curve that connects
$\phi_1^{-1}$
and
$\phi_2$.
Conversely, given a path 
$c: [0,1 ]\to G$
connecting
$\phi_1^{-1}$
and
$\phi_2$,
one can view it as a union of two curves that join each other
at their common endpoint. One of these curves
can be taken for
$a_2$
and the group inverse of the other one for
$a_1$.
Appropriately parameterized these two new curves
form a system of paths belonging to
${\cal G} ({\cal C}_1,{\cal C}_2)$.
One easily sees that there is a one-to-one correspondence
between homotopy classes of elements of
${\cal G} ({\cal C})$
and
${\cal G}^\prime ({\cal C})$.

Given a connection 
$\nabla\in {\cal L} ({\cal C})$
and a trivialization
$\tilde{P}$
of 
$P\to \Sigma$, 
one can associate to it a path
$a_{\nabla, \tilde{P}}$
belonging to
${\cal G}^\prime ({\cal C})$.
Namely let
$\varphi_1^{-1}\in {\cal C}_1^{-1}$
and
$\varphi_2\in {\cal C}_2$
be the holonomies of
$\nabla$
(with respect to the fixed trivialization
$\tilde{P}$) 
along, respectively, the closed paths
${\{ 0\times t\}}_{0\leq t\leq 1}$
and
${\{ 1\times t\}}_{0\leq t\leq 1}$.
Define a path 
$a_{\nabla, \tilde{P}}: [0,1]\to G$ 
by taking the holonomy of
$\nabla$
along the closed path
${\{ t\to s\times t\}}_{0\leq t\leq 1}$
as
$a_{\nabla, \tilde{P}} (s)$.
The homotopy class of 
$a_{\nabla, \tilde{P}}$
depends only on
$[\nabla ]$
and will be denoted by
$\{ \nabla \} $.
The corresponding connected component of
${\cal G}^\prime ({\cal C})$
will be denoted by
${\cal G}_{\{ \nabla \}  }^\prime ({\cal C})$.
We will show in the 
Section~\ref{subsect-pf-l-1-area-leq-distance-k-area-leq-coarse-distance-cylinder} 
that any homotopy class of paths from
${\cal G}^\prime ({\cal C})$
can be represented as
$\{ \nabla \} $
for some
$[\nabla ]$.

Now the theorem (for the case of a cylinder) 
becomes an immediate corollary of the following lemma.

\bigskip
\begin{lem1}
\label{lem-k-area-distance-cyl}
\begin{equation}
\label{lem-k-area-geq-coarse-distance-cylinder-pf}
\inf_{\scriptscriptstyle\nabla\in {\cal L}_{[\nabla ]} 
({\cal C})} 
\| L^\nabla\|
=
\inf_{\scriptscriptstyle a\in {\cal G}_{\{ \nabla \}  }^\prime ({\cal C})} 
{\hbox{\it coarse-length}}\ (a)
\end{equation}
\end{lem1}
\smallskip

\bigskip
\subsubsection{Proof of 
$\geq$
in
(\ref{lem-k-area-geq-coarse-distance-cylinder-pf})}
\label{subsect-pf-k-area-geq-coarse-distance-cylinder}

Let us take a connection 
$\nabla\in {\cal L} ({\cal C} )$.
Let us choose a global trivialization 
$\tilde{P}$
of 
$P\to \Sigma$
so that the holonomy of
$\nabla$
along any interval inside
$[0,1]\times 0$
is identity.
Given the trivialization, let
$\varphi_1^{-1}$
and
$\varphi_2$
be the holonomies of
$\nabla$
along the closed paths
${\{ t\to 0\times t\}}_{0\leq t\leq 1}$
and
${\{ t\to 1\times t\}}_{0\leq t\leq 1}$ 
respectively.
Then, as it was described above, 
$\nabla$
defines a path 
$a_{\nabla, \tilde{P}}: [0,1]\to G$
connecting
$\varphi_1^{-1}$
and
$\varphi_2$
and belonging to
${\cal G}_{\{  \nabla \}  }^\prime ({\cal C})$.

In order to prove
(\ref{lem-k-area-geq-coarse-distance-cylinder-pf})
it is enough to prove
\begin{equation}
\label{eqn-curvature-1}
\|L^\nabla\|
\geq 
\max_s \|\frac{d a_{\nabla,\tilde{P}} }{ds} (s)\|.
\end{equation}

To prove
(\ref{eqn-curvature-1})
let us denote by 
$\partial /\partial s$
and
$\partial /\partial t$
the standard vector fields on
$\Sigma = [0,1]\times S^1$
and by
$X$
and
$Y$
their horizontal lifts.
Let
$X^s$
and
$Y^t$
be the vertical components of the flows of
$X$
and
$Y$
respectively.
Set
$\phi_{s,t} = Y^t X^s$
(this map is defined on a domain which depends on
$s$
and
$t$).
Consider a family 
$v_{s,t}$
of vector fields on the fiber:
\[
v_{s,t} = \frac{\partial \phi_{s,t}}{\partial s} = Y^t_\ast X .
\]
Then, since the holonomy of
$\nabla$
along any interval inside
$[0,1]\times 0$
is identity, one has
\[
v_{s, 1} (s, 1) = \frac{d a_{\nabla, \tilde{P}}}{ds} (s)
\]
and
\[
v_{s, 0} (s, 0) = 0.
\]

On the other hand,
\[
\frac{\partial v_{s,t}}{\partial t} = Y^t_\ast [X,Y]^{vert},
\]
where
$[X,Y]^{vert}$
is the vertical component of
$[X,Y]$.
Thus
\[
\frac{d a_{\nabla, \tilde{P}}}{ds} (s) = 
v_{s,0} (s,0) + \int_0^1 \frac{\partial v_{s,1}}{\partial t} (s,0) dt = 
\int_0^1
Y^t_\ast [X,Y]^{vert} (s,0) dt,
\]
and hence
\begin{equation}
\label{eqn-curvature-2}
\int_0^1 \| Y^t_\ast [X,Y]^{vert} (s,0)\| dt \geq
\|\frac{d a_{\nabla,\tilde{P}} }{ds} (s) \|.
\end{equation}
The vector field 
$[X,Y]^{vert}$,
viewed by means of the trivialization of the bundle as an element
of
${\mathfrak g}$,
is equal (up to sign) to the curvature
$L^\nabla (\partial/{\partial s},\partial/{\partial t})$.
Since the Finsler norm on
$G$
is 
{\it
bi-invariant
}
we get that
\begin{equation}
\label{eqn-curvature-3}
\| Y^t_\ast [X,Y]^{vert} (s,0)\| \leq 
\|L^\nabla\|.
\end{equation}
Combining 
(\ref{eqn-curvature-2})
and
(\ref{eqn-curvature-3})
we get
(\ref{eqn-curvature-1})
and thus prove  
(\ref{lem-k-area-geq-coarse-distance-cylinder-pf}).
\b
\smallskip

\bigskip
\subsubsection{Proof of 
$\leq$
in
(\ref{lem-k-area-geq-coarse-distance-cylinder-pf}).
}
\label{subsect-pf-l-1-area-leq-distance-k-area-leq-coarse-distance-cylinder}

Let us take a path
$a: [0,1]\to G$
belonging to
${\cal G}^\prime ({\cal C})$.
Reparameterize
$a$
so that it becomes constant near 0 and 1.
Let
$K = [0,1]\times [0,1]$.
The trivial bundle
$P\to\Sigma$
can be constructed by taking the bundle
$K\times F \to K$
and identifying
$(s,0)\times a(s) y$
with
$(s,1)\times y$
for all
$s\in [0,1]$,
$y\in G$.
(Here
$a(s) y$
denotes the action of the element
$a(s)\in G$
on
$y\in F$).

We are going to define a connection
$\nabla$
on
$K\times F \to K$
by defining the corresponding field of horizontal planes.
Take a monotone cut-off function
$\psi (t)$
on the segment
$[0,1]$
such that
$\psi (t) = 1$
when 
$t$
is near 
$0$, 
and
$\psi (t) = 0$
when 
$t$
is near 
$1$.

The derivative
$\displaystyle{\frac{da}{ds}} (s)$
is a vector tangent to
$G$
at the point
$a (s)\in G$.
It can be identified with an element
$a^\prime (s)\in {\mathfrak g}$
corresponding to the right-invariant vector field on
$G$
containing
$\displaystyle{\frac{da}{ds}} (s)$.
On the other hand, if one considers the action of the elements
$a(s)$
on
$F$ 
then the vector
$\displaystyle{\frac{da}{ds}} (s)$
can be viewed as a vector field on
$F$
whose value at a point 
$y$
will be denoted by
$A (s,y)$.

Our horizontal plane will be generated by two horizontal
vector fields
$X$
and
$Y$
which are the horizontal lifts of
$\partial /\partial s$
and
$\partial /\partial t$ 
respectively.
The vector fields are defined as follows: at a point
$(s,t)\times y\in K\times F$
our horizontal plane will be generated by the vectors
\[
X_{(s, t)\times y} =\langle 
1,0,\psi(t) A (s,y)
\rangle
\]
and
\[
Y_{(s,t)\times y}=\langle
0,1,0
\rangle.
\]

One easily checks that the plane field spanned by
$X$
and
$Y$
defines a 
\break
$G{\hbox{\rm -connection}}$
$\nabla$
on
$P\to\Sigma$
which is flat near the boundary and
whose holonomies over the ends of
$\Sigma$
belong to the conjugacy classes 
${\cal C}_1$
and
${\cal C}_2$.
Therefore
$\nabla\in {\cal L} ({\cal C})$
and clearly
$[a] = \{  \nabla \} $.
This shows in particular that any homotopy class of paths from
${\cal G}^\prime ({\cal C})$
can be realized as 
$\{  \nabla \} $
for some
$[\nabla ]$.

The commutator of
$X$
and
$Y$
is a vertical vector field which at a point
$(s,t)\times y\in K\times F$
looks as follows:
\[
[X,Y]_{(s,t)\times y} = 
\langle
0,0, \psi^\prime (t) A (s,y)
\rangle .
\]

The curvature 
$L^\nabla ({\partial}/{\partial s},{\partial}/{\partial t})$
at the point
$(s,t)$
is the element
$a^\prime (s)\in {\mathfrak g}$.
Therefore, since the Finsler norm on
$G$
is 
{\it
bi-invariant
}
one has
\[
\|
L^\nabla (s,t) 
\|
= |\psi^\prime (t)|\cdot \| a^\prime (s)\|
= |\psi^\prime (t)|\cdot \|\frac{da}{ds} (s)\|.
\]
Hence
\[
\|L^\nabla\| = \max_{(s,t)} \|L^\nabla (s,t)\|
\leq
\max_t {|\psi^\prime (t)|} \cdot 
\max_s \|\frac{da}{ds} (s) \|.
\]

Now fix a small positive number
$\varepsilon$.
We can choose the function
$\psi$
above so that
$|\psi^\prime (t)|\leq 1+\varepsilon$
for all
$t$. 
Therefore we get
\[
\|L^\nabla\|
\leq
(1+\varepsilon)\ {\hbox{\it coarse-length}}\, (a).
\]

Since 
$\varepsilon$
was taken arbitrarily we get the desired inequalities. 
This finishes the proof of
Lemma~\ref{lem-k-area-distance-cyl}
and the proof of 
Theorem~\ref{thm-k-area-distance}
in the case when
$\Sigma$
is a cylinder.
\b

\bigskip
\begin{rem}
\label{rem-rescaling-cyl}
{\rm
Suppose that one drops the normalization condition for the area of the cylinder
$\Sigma$ 
to be 1. Then one can repeat the proofs above for the cylinder
$[0,A]\times S^1$
of area
$A$
taking the rescaling into account.
As a result one would get the inequalities:
\[
A\cdot 
\inf_{\scriptscriptstyle\nabla\in {\cal L} ({\cal C})} \| L^\nabla\|
\geq
\inf_{\scriptscriptstyle a\in {\cal G}^\prime ({\cal C})} 
{\hbox{\it length}}\,(a).
\]
}
\end{rem}
\smallskip

\bigskip
\subsection{The case 
$l\neq 2$}
\label{subsect-pf-k-area-dist-l-neq-1}

\bigskip
\subsubsection{The case of a disc:
$l = 1$}
\label{subsubsect-pf-k-area-dist-disc}

Let
$\Sigma$
be a disc. 
Present it as a cylinder 
$Cyl$
one of whose boundary components is 
capped with a disc
$D$. 
Considering trivial flat connections on
$P\to\Sigma$
defined over
$D$,
gluing them with connections defined over
$Cyl$
to get a connection over the whole
$\Sigma$
and using
Lemma~\ref{lem-k-area-distance-cyl}
and
Remark~\ref{rem-rescaling-cyl}
one easily proves the theorem for the case
$l=1$.
\b
\smallskip

\bigskip
\subsubsection{The case
$l > 1$}
\label{subsubsect-pf-k-area-dist-l>1}

Let us cut
$\Sigma$
into cylinders 
$Cyl_1, \ldots, Cyl_l$
(with piecewise smooth boundaries) in the following way.
Each boundary component of
$\Sigma$
will be a boundary component of exactly 
one of these cylinders.
If
$l=1$
we view the disc
$\Sigma$
as a cylinder with one boundary component capped with a disc.
Otherwise, if
$l >1$,
we make the cuts in such a way that the boundaries of 
the cylinders are not smooth only at some two common points
$p,\bar{p}\in\Sigma$.
The union
of all such boundaries passing through 
$p$
and
$\bar{p}$
consists of
$l$
paths
$\theta_1,\ldots,\theta_l$
coming out in the counterclockwise order from
$p$
and connecting it with
$\bar{p}$.
The closed paths
$c_1 = \theta_1\circ \theta_{l}^{-1}, 
c_2=\theta_2\circ \theta_1^{-1},\ldots,
c_l =\theta_l\circ \theta_{l-1}^{-1}: S^1\to\Sigma$
are, respectively, the boundary components of the
cylinders 
$Cyl_1, \ldots, Cyl_l$
(the other boundary components of these cylinders are, respectively,
$T_1,\ldots, T_l$).
Without loss of generality we can assume that the areas of
of the cylinders
$Cyl_1, \ldots, Cyl_l$
are all equal to
$\displaystyle{\frac{1}{l}}$.

Consider the standard cylinder
$[0, 1]\times S^1$
with the coordinates
$(s,t)$,
$0\leq s\leq 1$,
$0\leq t\leq 1\, ({\hbox{\rm mod}}\, 1)$,
equipped with the area form
$\displaystyle{\frac{ds\wedge dt}{l}}$
so that the total area is
$1/l$.

For each
$i=1,\ldots,l$
let us fix a point 
$y_i\in T_i$,
and a map
${\Phi}_i : [0,1]\times S^1\to Cyl_i$
satisfying the following conditions:

\medskip
\noindent
$\bullet$
each
${\Phi}_i : [0,1]\times S^1\to Cyl_i$
extends the corresponding map
$\Phi_i: [1-\delta, 1]\times S^1 \to\Sigma$
from
Section~\ref{subsect-surface};

\smallskip
\noindent
$\bullet$
${\Phi}_i (0\times t) = {\Phi}_{i+1} (0\times \{1 - t\}) = c_i (t)$,
$i=1,\ldots , l-1$,
${\Phi}_l (0\times t) = {\Phi}_1 (0\times \{ 1-t\})$;

\smallskip
\noindent
$\bullet$
${\Phi}_i (1\times 0) = y_i$,
${\Phi}_i (0\times 0) = p$;
${\Phi}_i (0\times 1/2) ={\bar p}$;

\smallskip
\noindent
$\bullet$
the map
${\Phi}_i$
is an
{\it
area-preserving
}
diffeomorphism over the pre-image of
$Cyl_i\setminus \{ p,{\bar p}\}$
and a homeomorphism over the closed cylinder.
\medskip

To a connection 
$\nabla$
on
$P = \Sigma\times G\to F$
with a chosen trivialization
$\tilde{P}$
one can associate an
$l{\hbox{\rm -system}}$
of paths 
${a}_{\nabla, \tilde{P}}$
in
$G$.
Namely, suppose 
$\nabla\in {\cal L} ({\cal C})$.
Set
$a_i (s)$
to be the parallel transport of
$\nabla$
along the path
$t\to \Phi_i (s\times t)$,
$i=1,\ldots,l$.
One can check that the paths
$a_i (s)$,
$0\leq s\leq 1$,
$i=1,\ldots,l$,
together form an
$l{\hbox{\rm -system}}$
of paths
${a}_{\nabla,\tilde{P}} \in {\cal G} ({\cal C})$.

As we mentioned in
Section~\ref{subsect-systems-of-paths}
the space 
${\cal G} ({\cal C})$ 
might have more than one connected component. 
One can check that if instead of 
$\tilde{P}$
we take another trivialization then the new system of
paths from
${\cal G} ({\cal C})$
is homotopic to the old one.
Thus the homotopy class of
${a}_{\nabla, \tilde{P}}$
depends only on the homotopy class
$[\nabla]$
and will be denoted by 
$[[\nabla]] = [a] ([\nabla])$. 

Gluing connections on
$P\to G$
defined over the cylinders
$Cyl_1,\ldots , Cyl_l$
into a connection defined over the whole
$\Sigma$
one can show that any homotopy class
$[a]$
can be represented as
$[[\nabla]] = [a] ([\nabla])$
for some
$[\nabla ]$.
Then using the same gluing technique together with
Lemma~\ref{lem-k-area-distance-cyl}
and
Remark~\ref{rem-rescaling-cyl}
one easily proves the theorem for the case
$l>1$.

This finishes the proof of  
Theorem~\ref{thm-k-area-distance}.
\b
\smallskip

\vfil
\eject

\bibliographystyle{alpha}

\end{document}